\documentclass[11pt,twoside]{article}
\pdfoutput=1 
\usepackage{epsfig,amsfonts,color}
\usepackage{amsmath}
\usepackage{mathtools}

\usepackage{amssymb, palatino, geometry,url}
\usepackage{amsthm}
\usepackage{algorithmic}
\usepackage{booktabs}
\usepackage{siunitx}
\usepackage{xcolor}
\usepackage{url}
\usepackage[most]{tcolorbox}
\usepackage[T1]{fontenc}
\newtcolorbox{highlighted}{colback=cyan,coltext=black,breakable}

\usepackage[noresetcount,lined,boxed]{algorithm2e} 
\newcommand{\mynote}[1]{{ \color{black}{#1}}}
\usepackage[colorlinks=true,linkcolor=blue,citecolor=blue,urlcolor=blue]{hyperref}
\usepackage{subcaption}
\usepackage{cleveref}
\usepackage{color}
\usepackage{fancyhdr}
\theoremstyle{definition}

\newcommand{\be}{\begin{equation}}
\newcommand{\ee}{\end{equation}}
\newcommand{\bea}{\begin{eqnarray}}
\newcommand{\eea}{\end{eqnarray}}

\newcommand{\bvec}{\left(\begin{array}{c}}
\newcommand{\evec}{\end{array}\right)}
\newcommand{\bsub}{\begin{subequations}}
\newcommand{\esub}{\end{subequations}}

\usepackage[thinc]{esdiff}

\usepackage{lineno}
\usepackage{tabularx}
\usepackage{multirow}
\usepackage{adjustbox}
\usepackage[outdir=./Output]{epstopdf}

\begin{document}

\title{Solution of Large-Scale Supply Chain Models using\\ Graph Sampling \& Coarsening}

\author{Jiaze Ma and Victor M. Zavala\thanks{Corresponding Author: victor.zavala@wisc.edu}\\
 {\small Department of Chemical and Biological Engineering}\\
 {\small University of Wisconsin-Madison, 1415 Engineering Dr, Madison, WI 53706, USA}}
\date{}
\maketitle

\begin{abstract}
We present a graph sampling and coarsening scheme (gSC) for computing lower and upper bounds for large-scale supply chain models. An edge sampling scheme is used to build a low-complexity problem that is used to finding an approximate (but feasible) solution for the original model and to compute a lower bound (for a maximization problem). This scheme is similar in spirit to the so-called sample average approximation scheme, which is widely used for the solution of stochastic programs. A graph coarsening (aggregation) scheme is used to compute an upper bound and to estimate the optimality gap of the approximate solution.  The coarsening scheme uses node sampling to select a small set of support nodes that are used to guide node/edge aggregation and we show that the coarsened model provides a relaxation of the original model and a valid upper bound. We provide numerical evidence that gSC can yield significant improvements in solution time and memory usage over state-of-the-art solvers. Specifically, we study a  supply chain design model (a mixed-integer linear program) that contains over 38 million variables and show that gSC finds a solution with an optimality gap of $<0.5$\% in less than 22 minutes.  
\end{abstract}


\section{Introduction}

The rapid expansion of manufacturing and distribution systems drive the development of large-scale supply chain models  \cite{villa2002emerging,garcia2015supply}.  Despite advances in computing power and algorithms, the complexity of supply chain models continuously defies the capabilities of off-the-shelf tools and motivates the development of decomposition and approximation schemes \cite{grossmann2012advances}. 

Decomposition approaches take advantage of the underlying partially-connected nature of supply chain problems and, as the name suggests, aim to solve the original problem by decomposing it in smaller subproblems and by coordinating the solutions of such subproblems. Well-known decomposition paradigms include bilevel, Benders, and Lagrangian  decomposition. Bilevel decomposition has mainly focused on multi-period supply chain models and aims to decompose the problem over time  \cite{bok2000supply,iyer1998bilevel,dogan2006decomposition}. Benders decomposition has been primarily used for solving stochastic supply chain models (by decomposing over scenarios)  \cite{pishvaee2014accelerated,oliveira2014accelerating,uster2007benders,santoso2005stochastic}. 
Lagrangian decomposition is a flexible paradigm that has been used for solving stochastic supply chain problems by decomposing over the scenario space and for solving deterministic supply chain problems by decomposing them over space and time dimensions  \cite{jackson2003temporal,oliveira2013lagrangean,terrazas2011temporal,terrazas2011multiscale,terrazas2011temporal,sousa2011global,van2001lagrangean}. These decomposition approaches are capable of finding solutions for large-scale supply chain models with tight optimality gaps  \cite{jackson2003temporal,oliveira2013lagrangean,sousa2011global}. 

Approximation is another paradigm that can be used to find solutions for large-scale supply chain problems; under this paradigm, one aims to solve an approximate/surrogate version of the original problem and to estimate the accuracy of the approximate solution. Aggregation (coarsening) is a powerful approach that is used to compute approximate solutions for optimization problems. At a fundamental level, coarsening omits local details of the problem while preserving the global structure \cite{barmann2015solving}. For instance, one can partition the elements of an optimization model (e.g., nodes in an underlying graph or network) according to attributes. The aggregation process then aims to aggregate/cluster elements with similar attributes by minimizing a distance measure \cite{chen2017large}. This distance can be defined, for instance, as the geographical distance between nodes in a graph \cite{barmann2015solving} or demands of products \cite{wu2021predictive}. Product aggregation has been used in inventory management problems \cite{chen2017large} and production planning problems \cite{zipkin1982exact}.  Network aggregation has been used in network flow problems \cite{zipkin1982transportation}, vehicle routing problems \cite{clautiaux2017iterative}, and network design problems \cite{barmann2015solving}). Aggregation has also been used to tackle general linear programs by aggregating constraints \cite{shetty1987solving} and/or variables \cite{zipkin1980variable}. This approach is also known as row aggregation and column aggregation \cite{litvinchev2013aggregation}. These aggregation schemes have also been used for computing bounds and for estimating the accuracy of approximate solutions (e.g., in the form of optimality gaps)  \cite{zipkin1980nodes}. 

A well-known and powerful approximation paradigm for the solution of stochastic programs (such as stochastic supply chain models) is sample average approximation (SAA) \cite{verweij2003sample}. This is a random sampling schemes that  deals with complexity associated with a large (or even infinite) number of scenarios by building an approximate model that uses a smaller scenario set \cite{kleywegt2002sample}. SAA has been widely used for solving large-scale stochastic supply chain models  in a wide range of applications  \cite{santoso2005stochastic,schutz2009supply,bidhandi2017accelerated,li2018sample}. A powerful feature of SAA is that it enables the computation of statistical lower and upper bounds of the true solution, which can be used to estimate the optimality gap. Diverse studies have shown that high-quality solutions can be obtained with this approach \cite{santoso2005stochastic,schutz2009supply}. 

In this work, we propose an approximation scheme based on graph sampling and coarsening (gSC) for computing lower and upper bounds for large-scale supply chain models and to estimate their solution accuracy. Our approach targets a deterministic design model for multi-product supply chains; this model uses a graph representation of the supply chain that captures spatial and product interdependencies between participants that arise from product transportation and transportation. The proposed scheme is motivated by the observation that the size of the underlying  graph of such problems scales (in the worst case) as O($|\mathcal{P}|\cdot |\mathcal{N}|^{2}$), where  $|\mathcal{N}|$ is the number of nodes and $|\mathcal{P}|$ is the number of products. This gives problems with a huge number of variables and constraints and hinders the ability of off-the-shelf solvers to tackle real applications. In addition to dimensionality, large-scale supply chain models are difficult to solve due to inherent solution degeneracy and tight algebraic coupling. The proposed gSC scheme is inspired by SAA in that is uses random edge sampling to build an approximate graph that uses a smaller (but representative) edge set. We show that this {\em graph-sampling} approach provides an approximate (but feasible) solution and provides a lower bound (for a maximization problem). To the best of our knowledge, this graph-sampling approach has not been considered in the literature. We also derive a {\em graph-coarsening} scheme for computing an upper bound and to estimate the optimality gap of the approximate solution; here, we borrow coarsening schemes previously used in the literature but adapt this to our context and extend these by incorporating edge aggregation capabilities.  Edge aggregation allows us to construct compact and tractable models and, to the best of our knowledge, this approach has not been considered in the literature. 
\\

We demonstrate the effectiveness of the proposed gSC scheme using a challenging agricultural waste management problem \cite{sampat2019coordinated}. This is a supply chain design problem that explores optimal transportation and transformation pathways of dairy waste to obtain diverse value-added products. The supply chain model contains 1,372 nodes, 1.8 million transport edges, 12 processing technologies, and 20 products and yields a mixed-integer linear program (MILP) with 38 million variables and half a million of constraints. Using the proposed approach, we find an approximate solution with an optimality gap of 0.5$\%$ in 22 minute. In contrast, an off-the-shelf solver (Gurobi) struggles to find solutions in reasonable times due to memory constraints. Our results suggest that the proposed scheme can be adapted to tackle other supply chain models. 

\section{Supply Chain Model}

We aim to develop algorithms for finding approximate solutions for large-scale multi-product supply chain models. The supply chain modeling framework considered here comprises a set of nodes $\mathcal{N}$ (spatial locations), suppliers $\mathcal{S}$, consumers $\mathcal{D}$, products $\mathcal{P}$, technologies $\mathcal{T}$, and transportation edges (links) $\mathcal{L}$ \cite{tominac2021economic}. The supply chain connects participants (suppliers, consumers, technologies, and transport services) through a complex graph. Specifically, the transportation participants offer product movement services across nodes via transportation edges, which induces spatial coupling. Technology participants offer processing services that convert raw materials into other value-added products (intermediate and final products),  which induces coupling across products.   

Each supplier $i\in\mathcal{S}$ has a supply flow $s_i\in\mathbb{R}_{+}$, a supplied product type $p(i)\in\mathcal{P}$, maximum supplying capacity $\bar{s}_i\in\mathbb{R}_+$, location $n(i)\in\mathcal{N}$, and cost $\alpha^s_i\in\mathbb{R}_+$. Each consumer $j\in\mathcal{D}$ has a demand flow $d_j\in\mathbb{R}_+$, requested product type $p(j)\in\mathcal{P}$, maximum demanding capacity $\bar{d}_j\in\mathbb{R}_+$, location of consumer $n(j)\in\mathcal{N}$, and cost $\alpha_j^d\in\mathbb{R}_+$. We use attributes to define the nested sets $\mathcal{S}_{n,p}\subseteq \mathcal{S}_n\subseteq\mathcal{S}$ with $\mathcal{S}_{n}:=\{i\in\mathcal{S} \mid n(i)=n \}$ (i.e., all suppliers attached to node $n$) and $\mathcal{S}_{n,p}:=\{i\in\mathcal{S} \mid n(i)=n, p(i)=p\}$ (i.e., all suppliers attached to node $n$ and providing product $p$ ). We follow a similar reasoning to define the nested sets $\mathcal{D}_{n,p}\subseteq\mathcal{D}_n\subseteq\mathcal{D}$.

Each technology $t\in \mathcal{T}$ has a set of yield factors $\gamma_{t,p}\in\mathbb{R}$, location $n(t)\in \mathcal{N}$, reference product $p(t)\in\mathcal{P}$, processing capacity $\bar{\xi}_{t}\in\mathbb{R}_+$, operating cost $\alpha^\xi_t\in \mathbb{R}_+$, number of facilities installed $y_t\in \mathbb{Z}_+$, number of facilities available $\bar{y}_{t}\in\mathbb{Z}_+$, installation cost $\alpha^y_t\in \mathbb{R}_+$, and node $n(t)\in \mathcal{N}$. Yield factors denote the units of product $p$ consumed/generated per unit of reference product $p(t)$ consumed/generated in the technology. A positive yield factor indicates product $p$ is generated in the technology $t$. A negative yield factor denotes the product $p$ is consumed. When the yield factor is zero, then product $p$ is neither produced nor consumed. The total amount of reference product $p(t)$ processed at technology $t$ is denoted as $\xi_{t}\in \mathbb{R}_+$. We denote the subset of technologies located in node $n$ as $\mathcal{T}_{n}\subseteq \mathcal{T}$, with $\mathcal{T}_{n}:=\{t\in\mathcal{T} \mid n(t)=n \}$. 

Each transportation edge $\ell \in\mathcal{L}$ has an associated flow $f_\ell\in\mathbb{R}_{+}$, product type transported $p(\ell)\in\mathcal{P}$, maximum capacity $\bar{f}_\ell\in \mathbb{R}_+$, and cost $\alpha^f_\ell\in \mathbb{R}_+$, sending (origin) node $n_s(\ell)\in\mathcal{N}$, and receiving (destination) node $n_r(\ell)\in\mathcal{N}$. We define the set of all flows entering node $n\in\mathcal{N}$ as $\mathcal{L}^{in}_n:=\{\ell\in\mathcal{L} \mid n_r(\ell)=n\}$ and we define the set of all edges leaving node $n\in\mathcal{N}$ as $\mathcal{L}^{out}_n:=\{\ell\in \mathcal{L} \mid n_s(\ell)=n\}$. We also define the nested subsets for entering edges ${\mathcal{L}^{in}_{n,p}}\subseteq \mathcal{L}^{in}_n\subseteq{\mathcal{L}}$ where $\mathcal{L}^{in}_{n,p}:=\{\ell\in\mathcal{L} \mid n_r(\ell)=n, p(\ell)=p\}$.  We use similar definitions to construct subsets $\mathcal{L}^{out}_{n,p}\subseteq \mathcal{L}^{out}_{n}\subseteq{\mathcal{L}}$.

\begin{figure}[!ht]
\begin{center}
\includegraphics[scale=0.3]{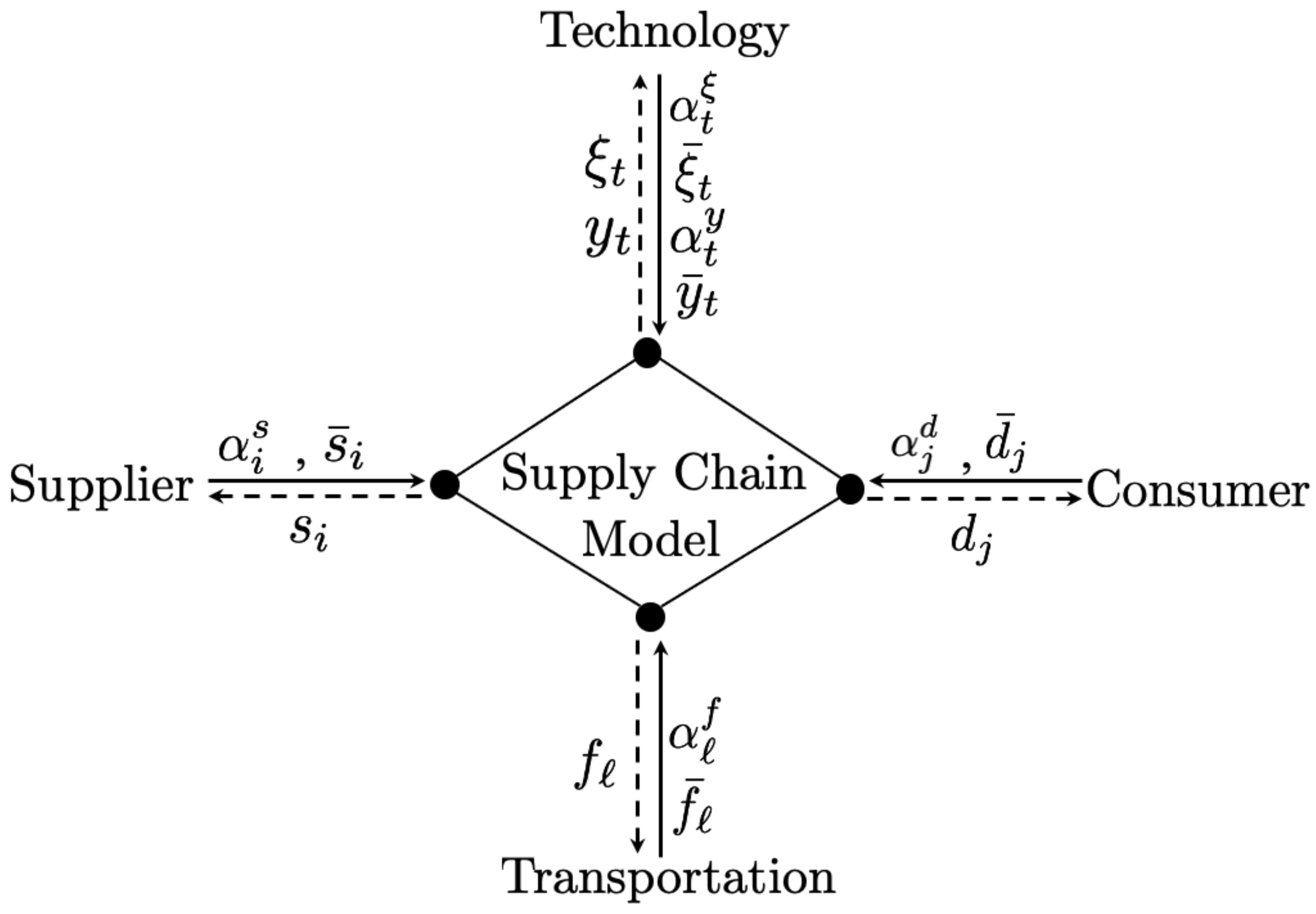}\caption{\mynote{Schematic of the supply chain model. Suppliers, consumers, and technology/transportation participants provide data consisting of cost bids $(\alpha^d, \alpha^s, \alpha^f,\alpha^\xi,\alpha^y_t)$ and capacities $(\bar{d}, \bar{s}, \bar{f},\bar{\xi},\bar{y})$. The supply chain model finds optimal allocations $(d, s, f,\xi,y)$ that maximize the total welfare.}}
\label{Fig.schematic_sc}
\end{center}
\end{figure}

The supply chain design model is given by: 
\begin{align}
\label{obj:social_welfare} \max_{(s,d,f,\xi,y)} & \; \sum_{j \in \mathcal{D}}\alpha^d_jd_j-\sum_{i \in \mathcal{S}}\alpha^s_is_i-\sum_{\ell \in \mathcal{L}}\alpha^f_\ell f_\ell-\sum_{t\in \mathcal{T}}\alpha^\xi_t\xi_t -\sum_{t\in \mathcal{T}}\alpha^y_t y_t\\ 
\label{eq:balancesp}\textrm{s.t.}& \left(\sum_{i\in \mathcal{S}_{n,p}}s_i  + \sum_{{\ell} \in \mathcal{L}^{in}_{n, p}}f_\ell \right) -\left(\sum_{j \in \mathcal{D}_{n, p}}d_j+\sum_{\ell \in \mathcal{L}^{out}_{n,p}}f_{\ell}\right) + \sum_{t\in \mathcal{T}_n}\gamma_{t,p}\, \xi_{t} =0,\;  (n,p)\in\mathcal{N}\times\mathcal{P}\; \\ 
&\; 0\leq d_j\leq \bar{d}_j,\; j\in \mathcal{D}\label{eq:cap1}\\
&\; 0\leq s_i\leq \bar{s}_i,\; i\in \mathcal{S}\label{eq:cap2}\\
&\; 0\leq f_\ell\leq \bar{f}_\ell,\; \ell\in \mathcal{L}\label{eq:cap3}\\
&\; 0\leq y_{t}\leq \bar{y}_{t},\; t\in \mathcal{T}.\label{eq:cap4}\\
&\; 0\leq \xi_t\leq y_t\bar{\xi}_t,\; t\in \mathcal{T}.\label{eq:cap5}
\end{align}
We refer to the objective function of this problem as the total welfare (denoted as $\varphi$) and we denote its optimal value as $\varphi^*$. As shown in Figure \ref{Fig.schematic_sc}, all participants submit data, consisting of costs $(\alpha^d, \alpha^s, \alpha^f,\alpha^\xi,\alpha^y_t)$ and capacity limits $(\bar{d}, \bar{s}, \bar{f},\bar{\xi},\bar{y})$. The supply chain model identifies optimal allocations $x=(d,s,f,\xi,y)$ that maximize the total welfare and that satisfy the product balance constraints \eqref{eq:balancesp}, capacity limits  \eqref{eq:cap1}-\eqref{eq:cap5}. Maximizing the total welfare \eqref{obj:social_welfare} aims to maximize the value fo the demand served while minimizing supply, technology, and transport costs. In the balance constraint \eqref{eq:balancesp}, the left parenthesis denotes the total input flow of product $p$ into node $n$, including the total supply flows $\sum_{i\in \mathcal{S}_{n,p}}s_i$, and inlet flows from other nodes $\sum_{{\ell} \in\mathcal{L}^{in}_{n, p}}f_\ell$. The right parenthesis denotes the total output flow of product $p$ from node $n$, including the total demand flows $\sum_{j\in \mathcal{D}_{n,p}}d_i$, and the outlet flows to other nodes $\sum_{{\ell} \in\mathcal{L}^{out}_{n, p}}f_\ell$. The third term accounts for the generation/consumption of product $p$ in all technologies $t$ located at node $n$. Constrains \eqref{eq:cap5} impose capacity bounds for each technology; the processing capacity of a technology is zero if the facility is not built ($y_t=0$), otherwise ($y_t>0$),  the total capacity is given by the per facility capacity ($\bar{\xi}_t$) and the number of facilities built ($y_t$). The feasible region defined by the constraints of the supply chain model is denoted as $\mathcal{F}$ and we note that this set contains the trivial allocation $(0,0,0,0,0)$ and we note that the supply chain problem is an MILP. 
\\

The supply chain model under study offers a compact representation of complex product transport and transformation tasks that take place in the system. This model also has interesting and important economic properties; specifically, one can show that supply chain design problem is equivalent to the design of a coordinated market and this enables the discovery of product prices and revenue streams allocated to participants \cite{tominac2021economic}. 

\section{Graph Sampling and Coarsening Scheme}

The supply chain model embeds a graph that scales (in the worst case) as O($|\mathcal{P}|\cdot |\mathcal{N}|^{2}$). This complexity arises because, in many applications of interest, it is possible to move products across any pair of nodes in the supply chain.  This complexity makes models of interest intractable to off-the-shelf solvers as this affects the number of decision variables, induces solution degeneracy, and induces  algebraic coupling (which affects linear algebra operations taking place inside solvers).  As such, it might be impossible to solve the supply chain model and compute the optimal welfare $\varphi^*$. The proposed approximation scheme (gSC) uses a graph sampling scheme that aims to obtain an approximate  solution and a valid lower bound for the optimal welfare (the supply chain model is a maximization problem) and uses a coarsening scheme to obtain a valid upper bound that helps estimate the accuracy (optimality gap) of the approximate solution. 

\subsection{Graph Sampling}

\begin{figure}[!ht]
\includegraphics[width=6in]{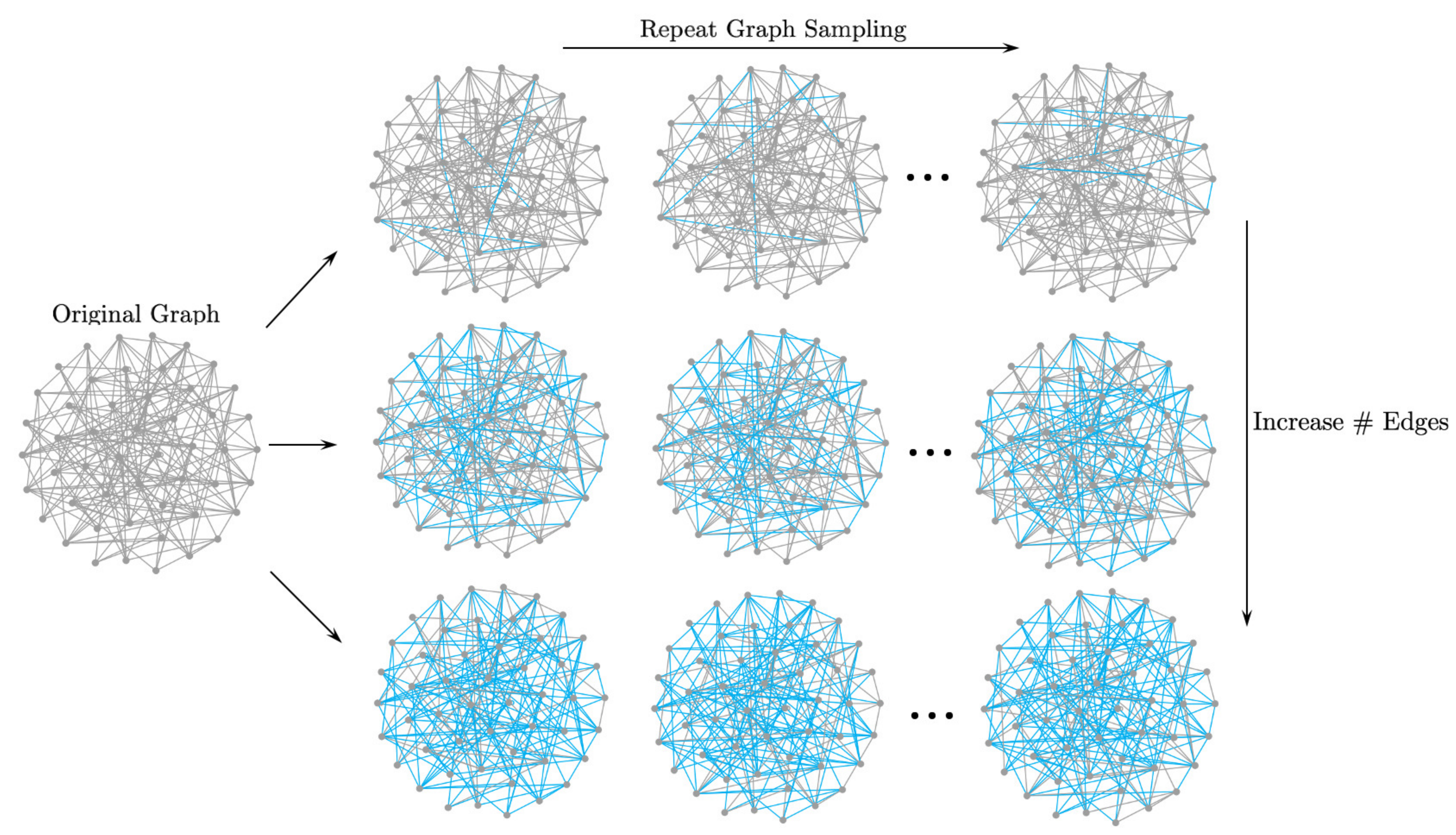}\caption{\mynote{Illustration of graph-sampling process. A small set of edges (in blue) are sampled and chosen to be active. Edges that are not sampled (in grey) are ignored. For a given sample size, the process is repeated to obtain a set of candidate solutions. Increasing the sample size tends to mitigate edge variance and associated variance in the solution.}}
\label{Fig.network_sampling}
\end{figure}

Inspired by SAA, we propose a graph-sampling approach that generates an approximate (smaller) graph of the supply chain. Here, we randomly sample a set of active edges $\mathcal{L}_a$ from the entire edge set $\mathcal{L}$. The edges that are not sampled are discarded and this equivalent to assume that these have zero flow. Edge sampling thus  restricts the set of available flows and thus the feasible region, yielding a suboptimal solution and providing a lower bound. Figure \ref{Fig.network_sampling} illustrates this scheme.  \\

\begin{figure}[!htp]
\begin{center}
\includegraphics[width=6in]{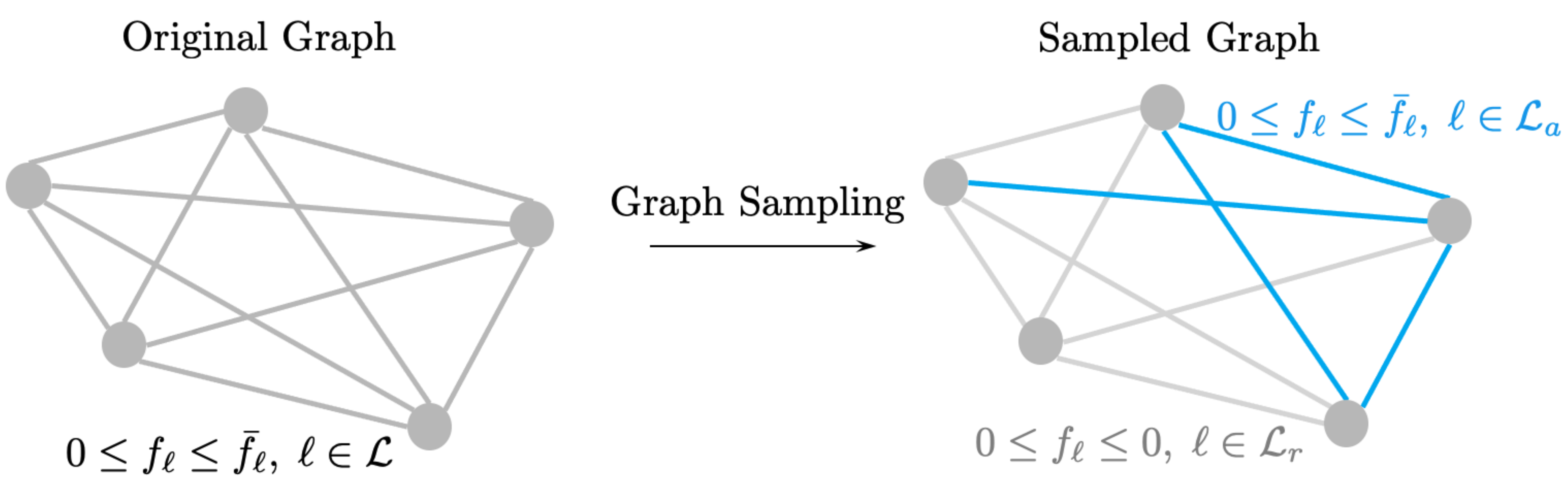}\caption{\mynote{\textcolor{black}{Illustration of graph sampling with notation. On the left is the original graph and $\mathcal{L}$ denotes the entire edge set.  On the right, $\mathcal{L}_a$ is the set of active edges (in blue) and $\mathcal{L}_r$ is the set of inactive edges (in grey),  there is no flow going through the inactive edges.}}}
\label{Fig.simple_sample}
\end{center}
\end{figure}

\textcolor{black}{Figure.  \ref{Fig.simple_sample} illustrates the graph sampling process.}The edge sampling procedure induces the set partition $\mathcal{L}=\mathcal{L}_a\cup\mathcal{L}_r$ with $\mathcal{L}_a\cap\mathcal{L}_r=\emptyset$ and express the flow limit constraints as:
\begin{subequations}
\begin{align}
& 0\leq f_{\ell}\leq \bar{f}_{\ell},\; \ell\in \mathcal{L}_a\\
& 0\leq f_{\ell}\leq 0,\; \ell\in \mathcal{L}_r.
\end{align}
\end{subequations}

Using the edge partition induced via sampling, we construct the approximate supply chain model:
\begin{subequations}
\begin{align}
\max_{(s,d,f,\xi,y)} & \; \sum_{j \in \mathcal{D}}\alpha^d_jd_j-\sum_{i \in \mathcal{S}}\alpha^s_is_i-\sum_{\ell \in \mathcal{L}}\alpha^f_\ell f_\ell-\sum_{t\in \mathcal{T}}\alpha^\xi_t\xi_t -\sum_{t\in \mathcal{T}}\alpha^y_t y_t\\ 
\textrm{s.t.}& \left(\sum_{i\in \mathcal{S}_{n,p}}s_i  + \sum_{{\ell} \in \mathcal{L}^{in}_{n, p}}f_\ell \right) -\left(\sum_{j \in \mathcal{D}_{n, p}}d_j+\sum_{\ell \in \mathcal{L}^{out}_{n,p}}f_{\ell}\right) + \sum_{t\in \mathcal{T}_n}\gamma_{t,p}\, \xi_{t} =0,\;  (n,p)\in\mathcal{N}\times\mathcal{P}\; \\ 
&\; 0\leq d_j\leq \bar{d}_j,\; j\in \mathcal{D}\\
&\; 0\leq s_i\leq \bar{s}_i,\; i\in \mathcal{S}\\
&\; 0\leq f_{\ell}\leq \bar{f}_{\ell},\; \ell\in \mathcal{L}_a\\
&\; 0\leq f_{\ell}\leq 0,\; \ell\in \mathcal{L}_r\\
&\; 0\leq y_{t}\leq \bar{y}_{t},\; t\in \mathcal{T}\\
&\; 0\leq \xi_t\leq \bar{\xi}_t\cdot y_t,\; t\in \mathcal{T}
\end{align}
\end{subequations}

We denote the feasible region of this approximate supply chain model is denoted as $\underline{\mathcal{F}}$ and the associated optimal total welfare as $\underline{\varphi}^*$. We note that the feasible regions $\underline{\mathcal{F}}$ and $\mathcal{F}$ (for the original supply chain model) both contain the trivial solution allocation $(0,0,0,0,0)$. Moreover, since the approximate model forces $f_{\ell}=0,\; \ell \in\mathcal{L}_r$, we have that $\underline{\mathcal{F}}\subseteq \mathcal{F}$. Therefore, the approximate model provides a feasible solution and a valid lower bound for the optimal total welfare of the original supply chain model ($\underline{\varphi}^*\leq \varphi^*$). 
\\

Given a desired number of active edges $a=|\mathcal{L}_a|\in \mathbb{Z}_+$ we can repeat the sampling process for multiple samples $\mathcal{L}_a(\omega),\; \omega \in \Omega$ and solve the approximate supply chain model using such samples to obtain samples for the total welfare $\underline{\varphi}_\omega,\,\omega \in \Omega$. Each of these welfare samples satisfies $\underline{\varphi}_\omega\leq \varphi$ and we can compute statistics for the welfare such as the expected value $\mathbb{E}[\underline{\varphi}]$ and standard deviation $\mathbb{SD}[\underline{\varphi}]$. A large variance indicates that the current sample is small and thus sampled graphs do not reflect the characteristics of the full original graph (there is a high variability in the graph topology). To mitigate this variability, we can increase the sample size and repeat the process; increasing the sample size is expected to reduce variability in the graph structure and in the associated total welfare. A small variance of the solution indicates that either the sampled graphs are representative of the full graph or that the problem exhibits solution degeneracy (different sample graphs and associated solutions yield a similar total welfare ). 

\subsection{Graph Coarsening}

Coarsening (aggregation) is a classical and powerful technique that has been used for solving large-scale network flow problems. The basic idea behind these models is to replace the large-scale graph (network) with a smaller, approximate graph by aggregating nodes and edges. Solving the coarsened  problem reduces computational complexity but the challenge is to ensure that the aggregated problem is constructed in a way that it provides a relaxation of the original  model and a valid upper bound for the optimal welfare $\varphi^*$.  
\\

The proposed aggregation scheme adapts partitioning and aggregation schemes proposed in \cite{barmann2015solving,zipkin1980nodes} to the multi-product supply chain setting of interest. Unfortunately, we will see that node aggregation is not sufficient to obtain a computationally tractable approximation; As such, we also propose to use an edge aggregation scheme.  The combination of node and edge aggregation will provide a tractable approximation and a valid upper bound. 

\begin{figure}[!htb]
\begin{center}
\includegraphics[width=6in]{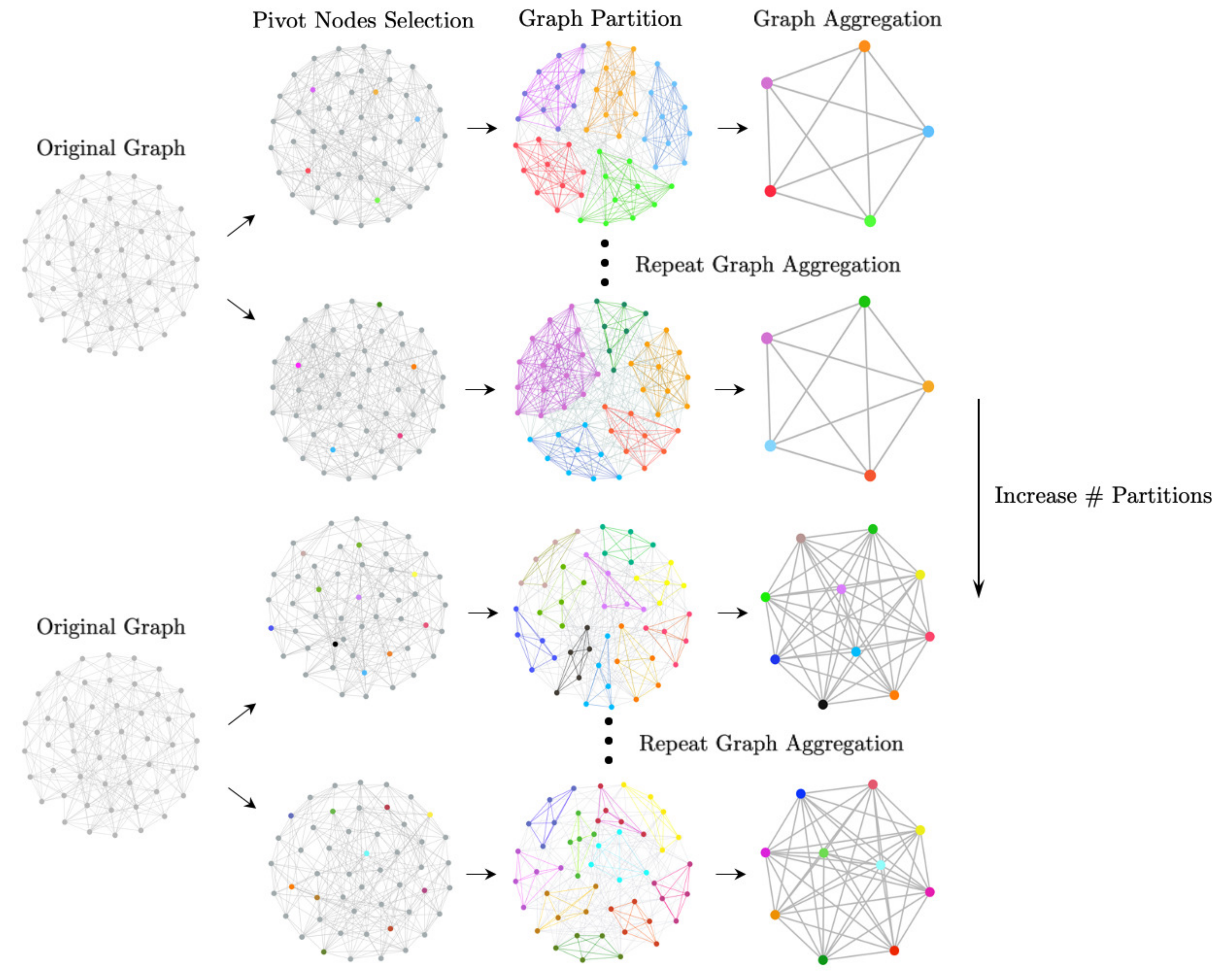}\caption{\mynote{Illustration of the graph aggregation process.  A small number of nodes are randomly chosen to be pivot nodes (colored nodes). The neighboring nodes (in grey) are aggregated with the corresponding pivot nodes to form graph partitions. The local (internal) edges in each partition are eliminated via node aggregation and the global edges connecting partitions are aggregated. The compact graph obtained is shown on the right. To increase the quality of the approximation, we increase the number of partitions (active nodes) and repeat the aggregation process.}}
\label{Fig.network_agg}
\end{center}
\end{figure}

\begin{figure}[!htp]
\begin{center}
\includegraphics[width=6in]{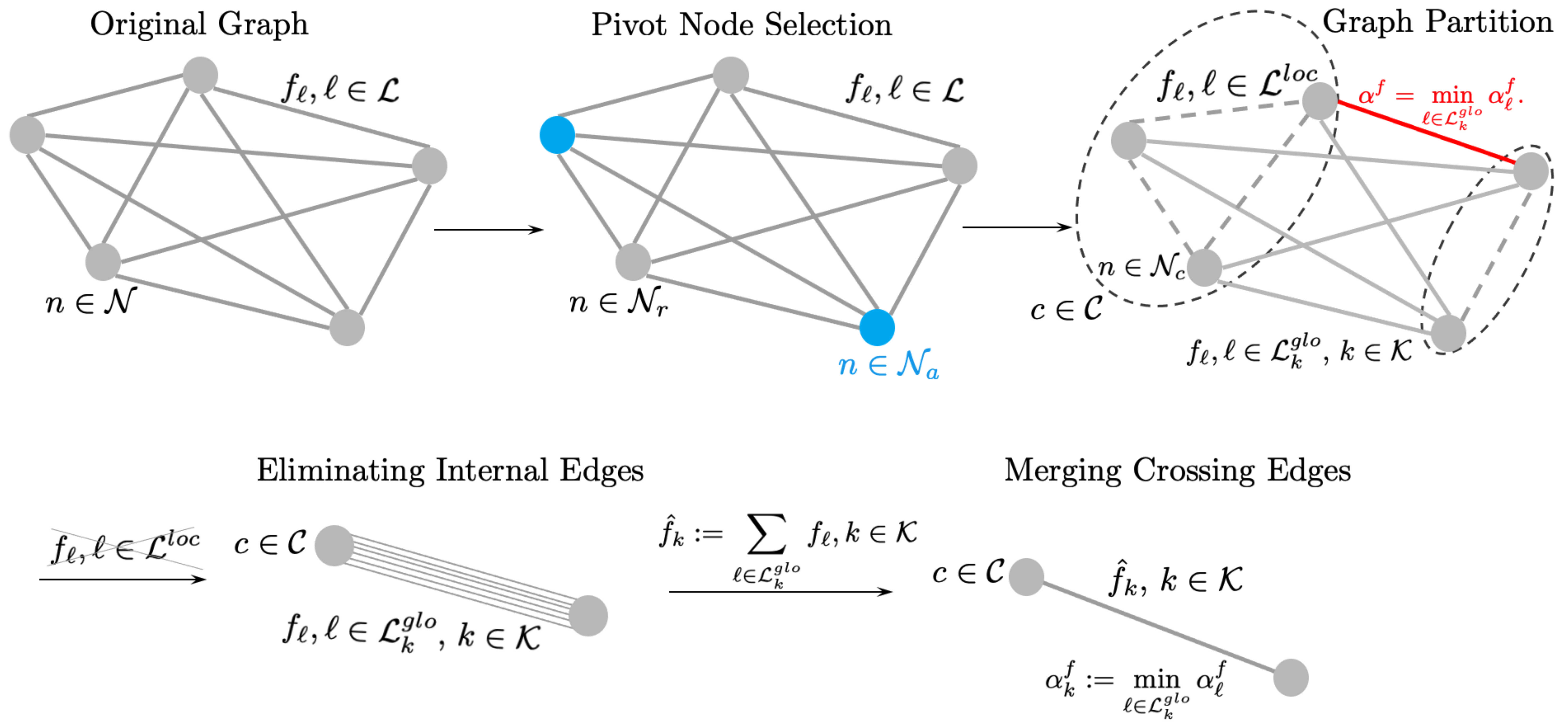}\caption{\mynote{\textcolor{black}{Illustration of graph coarsening with notations.  In the first plot, $\mathcal{N}$ is the entire node set and $\mathcal{L}$ is the entire edge set.  In the second plot,  the blue nodes are the pivot nodes and $\mathcal{N}_a$ is the pivot node set.  $\mathcal{N}_r$ is the remaining node set. In the third plot, $\mathcal{C}$ is the partition set,  $\mathcal{N}_c$ consists all the nodes within partition $c$,$c\in\mathcal{C}$,  $\mathcal{L}^{loc}$ is the set of local edges (dash lines) within each partition,  $\mathcal{L}^{glo}$ is the set of crossing edges (solid lines). Among the crossing edges $\ell \in \mathcal{L}^{glo}_k$, the one in red has the lowest flow cost $\min_{\ell\in\mathcal{L}_k^{glo}}$ $\alpha^f_{\ell}$. In the fourth plot, the internal edges $\ell\in\mathcal{L}_k^{loc}$ are eliminated; only the crossing edges $\ell\in\mathcal{L}_k^{glo}$ are preserved. In the fifth plot, the crossing edges are aggregated into a new edge, and the flow through this edge is named $\hat{f}_{k}, k\in\mathcal{K}$.  The flow cost of this edge is set to the lowest flow cost of the merged global edges, $\alpha^f_{k}:=\min_{\ell\in\mathcal{L}_k^{glo}} \alpha^f_\ell$.}}}
\label{Fig.simple_agg}
\end{center}
\end{figure}

Figures \ref{Fig.network_agg} and  \ref{Fig.simple_agg} illustrate the graph aggregation process. Given the supply chain graph, we begin by {\em randomly} selecting a small set of nodes $\mathcal{N}_a$ from the entire node set $\mathcal{N}$. The selected nodes are called \textit{pivot nodes} (or active nodes). Since the pivot nodes are selected at random, we also referred to them as sampled nodes. We define the complementary (remaining) set of nodes as $\mathcal{N}_r$ and we note that $\mathcal{N}_a\cup \mathcal{N}_r=\mathcal{N}$ and $\mathcal{N}_a\cap \mathcal{N}_r=\emptyset$.  The graph aggregation seeks to collapse all nodes into the pivot nodes; to do so, we begin by inducing a node set partition; here, we assign a pivot node $n\in \mathcal{N}_a$ to a node partition $c\in \mathcal{C}$ (with partition index $\mathcal{C}:=\{1,2,...,C\}$) and we assign nodes from $\mathcal{N}_r$ to such partitions. This assignment is cast as an optimization problem of the form: 
\begin{subequations}
\begin{align}
\label{obj:total_distance} \min_{z}\sum_{n\in \mathcal{N}_a}\sum_{n' \in \mathcal{N}_r} \theta_{n,n'}\cdot z_{n,n'}\\
\sum_{n\in \mathcal{N}_a}\ z_{n,n'}=1,\;n'\in\mathcal{N}_r \label{cons:logic}
\end{align}
\end{subequations}


The objective \eqref{obj:total_distance} is to minimize the total distance between the active nodes and their neighboring nodes that are aggregated with them. Distance is defined in terms of an attribute; here, we use the geographical distance between nodes (captured implicitly in the transportation costs $\alpha_\ell$). This is because, as we will see, the goal of node aggregation is to ignore local (short-distance) transport costs. The binary variable $z_{n,n'}\in \{0,1\}$ defines whether node $n'\in \mathcal{N}_r$ is aggregated with the pivot node $n\in \mathcal{N}_a$ and $\theta_{n,n'}\in \mathbb{R}_+$ denotes the distance between nodes $n\in \mathcal{N}_a$ and $n'\in \mathcal{N}_{r}$. The constraint \eqref{cons:logic} enforces that each node is only assigned to one partition (thus yielding non-overlapping partitions). For a given node $n\in \mathcal{N}$, we define its assigned partition as $c(n)\in \mathcal{C}$.
\\

The node partition induces a partition of the edge set as $\mathcal{L}=\mathcal{L}^{loc}\cup\mathcal{L}^{glo} = \mathcal{L}$ with $\mathcal{L}^{loc}\cap\mathcal{L}^{glo} = \emptyset$.  Here, $\mathcal{L}^{loc}$ is the set of local edges (internal to partitions) and $\mathcal{L}^{glo}$ is the set of global edges  (connect partitions). The edge partitions  are constructed as:
\begin{subequations}
\begin{align}
\mathcal{L}^{loc}&=\{\ell\in \mathcal{L}\,|\,c({n}_s(\ell))=c({n}_r(\ell))\}\\
\mathcal{L}^{glo}&=\{\ell\in \mathcal{L}\,|\,c({n}_s(\ell))\neq c({n}_r(\ell))\}. 
\end{align}
\end{subequations}
In other words, a local edge $\ell \in \mathcal{L}^{loc}$ has the same sending and receiving partitions $c({n}_s(\ell))\in \mathcal{C}$ and $c({n}_r(\ell))\in \mathcal{C}$. A global edge $\ell \in \mathcal{L}^{glo}$ has different sending and receiving node partitions $c({n}_s(\ell))\in \mathcal{C}$ and $c({n}_r(\ell))\in \mathcal{C}$. For convenience, we define the sending partition of an edge as $c_s(\ell):=c({n}_s(\ell))$ and the receiving partition of an edge as $c_r(\ell):=c({n}_r(\ell))$.
\\

As is typically done in aggregation schemes, we eliminate the local nodes and edges within each partition (and thus collapse the partition into its pivot node). This process results in a graph with $|\mathcal{C}|$ nodes and $|\mathcal{L}^{glo}|$ edges. This aggregation process has been applied to large-scale shortest path problems \cite{zipkin1980nodes} and network design problems \cite{barmann2015solving}. Unfortunately, this approach does not reduce the number of global edges and this hinders our ability to tackle problems of interest. Specifically, the model under study allows for transportation of multiple products between nodes (and thus a large number of global edges are typically encountered). 
\\

To further reduce the graph complexity, we propose an approach that also aggregates global edges. To do so, we partition the global edge set as $\mathcal{L}^{glo}=\cup_{k\in\mathcal{K}}\mathcal{L}_{k}^{glo}$ where $\mathcal{K}:=\{1,2,...,K\}$.  Here, all the edges in a given edge partition $\ell\in \mathcal{L}_k^{glo}$ have the same sending and receiving node partitions. To conduct global edge aggregation, we collapse a given edge partition $\mathcal{L}_k^{glo}$ into a single edge. Here, an aggregated edge $k\in \mathcal{K}$ has an associated flow $\hat{f}_k\in \mathbb{R}_+$ and has associated sending and receiving node partitions $c_s(k)\in \mathcal{C}$ and $c_r(k)\in \mathcal{C}$.  This approach reduces the number of edges to $|\mathcal{K}|$. 
\\

We now construct an approximate model that uses the aggregated graph representation containing $\mathcal{C}$ nodes and $\mathcal{K}$ edges.  The goal is to show that this aggregated formulation provides a relaxation of the original supply chain model and delivers a valid upper bound for the optimal welfare $\varphi^*$. 
\\

We first note that the flow constraint set of the original model can be expressed as:
\begin{subequations}
\begin{align}
& 0\leq f_\ell \leq \bar{f}_\ell,\; \ell\in\mathcal{L}^{loc}\\
& 0\leq f_\ell \leq \bar{f}_\ell,\; \ell \in \mathcal{L}_k^{glo}, \, k\in\mathcal{K}
\end{align}
\end{subequations}
We relax this constraint set by eliminating the local edges and by aggregating the flows of the global edges (via linear combination):
\begin{subequations}
\begin{align}
\label{eq:externalflow} 0\leq\sum_{\ell\in\mathcal{V}_k}f_l\leq \sum_{\ell\in\mathcal{V}_k}\bar{f}_{l},\ k\in\mathcal{K}. 
\end{align}
\end{subequations}
By defining the aggregated global flows and capacities as: 
\begin{subequations}
\begin{align}
\hat{f}_k& := \sum_{\ell\in\mathcal{L}_k^{glo}}f_\ell,
, k\in\mathcal{K}\\
\bar{\hat{f}}_k&:=\sum_{\ell\in\mathcal{L}_k^{glo}}\bar{f}_{\ell},\, k\in\mathcal{K}
\end{align}
\end{subequations}
we can express the flow constraints as:
\begin{subequations}
\begin{align}
\label{eq:externalflow} 0\leq\hat{f}_k\leq \bar{\hat{f}}_{k},\ k\in\mathcal{K}.
\end{align}
\end{subequations}
We now note that the nodal balance constraint set for the original model can be expressed as:
\begin{align}
\label{eq:newbalance_ub1}&\left(\sum_{i\in \mathcal{S}_{n,p}}s_i  + \sum_{{\ell} \in \mathcal{L}^{in}_{n, p}}f_\ell \right) -\left(\sum_{j \in \mathcal{D}_{n, p}}d_j+\sum_{\ell \in \mathcal{L}^{out}_{n,p}}f_{\ell}\right) + \sum_{t\in \mathcal{T}_n}\gamma_{t,p}\, \xi_{t} =0,\;  (n,p)\in\mathcal{N}_c\times\mathcal{P},\; c\in \mathcal{C}. 
\end{align}
We relax this constraint set by aggregating the nodal balances (via a linear combination) in each partition $c\in \mathcal{C}$ as:
\begin{align}
\label{eq:newbalance_ub2}&\left(\sum_{n\in {\mathcal{N}_c}}\sum_{i\in \mathcal{S}_{n,p}}s_i  +  \sum_{n\in {\mathcal{N}_c}}\sum_{{\ell} \in \mathcal{L}^{in}_{n,p}}f_\ell \right) -\left(\sum_{n\in {\mathcal{N}_c}}\sum_{j \in \mathcal{D}_{n,p}}d_j+\sum_{n\in {\mathcal{N}_c}}\sum_{\ell \in \mathcal{L}^{out}_{n,p}}f_{\ell}\right) + \sum_{n\in {\mathcal{N}_c}}\sum_{t\in \mathcal{T}_{n}}\gamma_{t,p}\, \xi_{t} =0,\;  (c,p)\in\mathcal{C}\times\mathcal{P}\;
\end{align} 
We now split the set of edges entering a node $n\in \mathcal{N}_c$ into the subsets $\mathcal{L}^{in}_{n,p}=\mathcal{L}^{in,loc}_{n,p} \cup \mathcal{L}^{in,glo}_{n,p}$. These subsets are defined as: 
\begin{align}
\mathcal{L}^{in,loc}_{n,p}&:=\{\ell \in \mathcal{L}\mid c_r(\ell)=c_s(\ell)\}\\
 \mathcal{L}^{in,glo}_{n,p}&:=\{\ell \in \mathcal{L}\mid c_r(\ell)\neq c_s(\ell)\}.
\end{align}
Similarly, the set of outlet edges is split as $\mathcal{L}^{out}_{n,p}=\mathcal{L}^{out,loc}_{n,p}\cup \mathcal{L}^{out,glo}_{n,p}$.  This allows us to express the aggregated balances as:
\begin{align}
\begin{split}
\label{eq:newbalance_ub3}&\left(\sum_{n\in {\mathcal{N}_c}}\sum_{i\in \mathcal{S}_{n,p}}s_i  +  \sum_{n\in {\mathcal{N}_c}}\sum_{{\ell} \in \mathcal{L}^{in,loc}_{n,p}}f_\ell +\sum_{n\in {\mathcal{N}_c}}\sum_{{\ell} \in \mathcal{L}^{in,glo}_{n,p}}f_\ell \right) -\left(\sum_{n\in {\mathcal{N}_c}}\sum_{j \in \mathcal{D}_{n,p}}d_j+\sum_{n\in {\mathcal{N}_c}}\sum_{\ell \in \mathcal{L}^{out,loc}_{n,p}}f_{\ell}+\sum_{n\in {\mathcal{N}_c}}\sum_{\ell \in \mathcal{L}^{out,glo}_{n,p}}f_{\ell} \right)\\
&+ \sum_{n\in {\mathcal{N}_c}}\sum_{t\in \mathcal{T}_{n}}\gamma_{t,p}\, \xi_{t} =0,\;  (c,p)\in\mathcal{C}\times\mathcal{P}\;
\end{split}
\end{align}

One can show that the aggregation of local flows entering the nodes in a partition equals to the aggregation of local flows leaving the nodes of a partition:
\begin{align}
\label{eq:newbalance_ub4}\sum_{n\in {\mathcal{N}_c}}\sum_{{\ell} \in \mathcal{L}^{in,loc}_{n,p}}f_\ell= \sum_{n\in {\mathcal{N}_c}}\sum_{{\ell} \in \mathcal{L}^{out,loc}_{n,p}}f_\ell,\;  (c,p)\in\mathcal{C}\times\mathcal{P}.
\end{align}
Using this fact, we can write the aggregated nodal balances: 
\begin{align}
\label{eq:newbalance_ub5}&\left(\sum_{n\in {\mathcal{N}_c}}\sum_{i\in \mathcal{S}_{n,p}}s_i+\sum_{n\in {\mathcal{N}_c}}\sum_{{\ell}\in\mathcal{L}^{in,glo}_{n,p}}f_\ell \right) -\left(\sum_{n\in {\mathcal{N}_c}}\sum_{j \in \mathcal{D}_{n,p}}d_j+\sum_{n\in {\mathcal{N}_c}}\sum_{\ell \in \mathcal{L}^{out,glo}_{n,p}}f_{\ell} \right)
+ \sum_{n\in {\mathcal{N}_c}}\sum_{t\in \mathcal{T}_{n}}\gamma_{t,p} \xi_{t} =0,\;  (c,p)\in\mathcal{C}\times\mathcal{P}.
\end{align}

We now note that the total number of suppliers, consumers, and technologies in each partition $c\in\mathcal{C}$ is simply the aggregation of such participants in the nodes associated with such partitions:
\begin{subequations}
\begin{align}
\label{eq:supplyagg}\sum_{n\in {\mathcal{N}_c}}\sum_{i\in\mathcal{S}_{n,p}} s_i&=\sum_{i\in\mathcal{S}_{c,p}} s_i,\;(c,p)\in\mathcal{C}\times\mathcal{P}\\  
\label{eq:demandagg}\sum_{n\in {\mathcal{N}_c}}\sum_{j\in\mathcal{D}_{n,p}} d_j&=\sum_{i\in\mathcal{D}_{c,p}} d_j,\;(c,p)\in\mathcal{C}\times\mathcal{P}\\  
\label{eq:processingagg}\sum_{n\in {\mathcal{N}_c}}\sum_{t\in \mathcal{T}_{n}}\gamma_{t,p}\xi_{t}&=\sum_{t\in \mathcal{T}_{c}}\gamma_{t,p}\xi_t,\;(c,p)\in\mathcal{C}\times\mathcal{P}
\end{align}
\end{subequations}
where:
\begin{subequations}
\begin{align}
\mathcal{S}_{c,p}&:=\{i\in\mathcal{S}\,|\, c(n(i))=c,\; p(i)=p\}\\
\mathcal{D}_{c,p}&:=\{j\in\mathcal{D}\,|\, c(n(j))=c,\; p(i)=p\}\\
\mathcal{T}_{c}&:=\{t\in\mathcal{T}\,|\, c(n(t))=c\}.
\end{align}
\end{subequations}
The aggregation of flows can now be expressed in terms of the aggregated global edges as:
\begin{subequations}
\begin{align}
\label{eq:flowsin}\sum_{n\in {\mathcal{N}_c}}\sum_{{\ell}\in\mathcal{L}^{in,glo}_{n,p}}f_\ell&=\sum_{{k} \in\mathcal{K}^{in}_{c, p}}\hat{f}_k,\; (c,p)\in\mathcal{C}\times\mathcal{P}\\
\label{eq:flowsout}\sum_{n\in {\mathcal{N}_c}}\sum_{{\ell}\in\mathcal{L}^{out,glo}_{n,p}}f_\ell&=\sum_{{k} \in\mathcal{K}^{out}_{c, p}}\hat{f}_k ,\;(c,p)\in\mathcal{C}\times\mathcal{P}
\end{align}
\end{subequations}
where:
\begin{subequations}
\begin{align}
\mathcal{K}^{in}_{c,p}&:=\{k\in \mathcal{K}\,|\,c_r(k)=c,\; p(k)=p\}\\
\mathcal{K}^{out}_{c,p}&:=\{k\in \mathcal{K}\,|\,c_s(k)=c,\; p(k)=p\}.
\end{align}
\end{subequations}
This yields the aggregated balance constraints:
\begin{align}
\label{eq:clusterbalancehat}&\left(\sum_{i\in\mathcal{S}_{c,p}}s_i + \sum_{{k} \in
\mathcal{K}^{in}_{c, p}}\hat{f}_k \right) -\left(\sum_{j\in\mathcal{D}_{c,p}}d_j+\sum_{k\in \mathcal{K}^{out}_{c,p}}\hat{f}_k\right) + \sum_{t\in \mathcal{T}_{c}}\gamma_{t,p}\xi_t =0,\; (c,p)\in\mathcal{C}\times\mathcal{P}.
\end{align}
\\

To aggregate the total welfare, we use the following sequence of observations:
\begin{subequations}
\begin{align}
\varphi
\label{eq:originalobj1_1} & =\sum_{j \in \mathcal{D}}\alpha^d_jd_j-\sum_{i \in \mathcal{S}}\alpha^s_is_i-\sum_{\ell\in\ \mathcal{L}}\alpha^{f}_\ell f_\ell-\sum_{t\in \mathcal{T}}\alpha^\xi_t\xi_t-\sum_{t\in \mathcal{T}}\alpha^y_t y_t\\
\label{eq:originalobj1_2}& =\sum_{j \in \mathcal{D}}\alpha^d_jd_j-\sum_{i \in \mathcal{S}}\alpha^s_is_i-\sum_{\ell\in\ \mathcal{L}^{loc}}\alpha^{f}_\ell f_\ell-\sum_{\ell\in\ \mathcal{L}^{glo}}\alpha^{f}_\ell f_\ell-\sum_{t\in \mathcal{T}}\alpha^\xi_t\xi_t - \sum_{t\in \mathcal{T}}\alpha^y_t y_t\\
\label{eq:originalobj1_3}& =\sum_{j \in \mathcal{D}}\alpha^d_jd_j-\sum_{i \in \mathcal{S}}\alpha^s_is_i-\sum_{\ell\in\ \mathcal{L}^{loc}}\alpha^{f}_\ell f_\ell-\sum_{k\in\ \mathcal{K}}\sum_{\ell\in\ \mathcal{L}_k^{glo}}\alpha^{f}_\ell f_\ell-\sum_{t\in \mathcal{T}}\alpha^\xi_t\xi_t - \sum_{t\in \mathcal{T}}\alpha^y_t y_t\\
\label{eq:originalobj3}&\leq \sum_{j \in \mathcal{D}}\alpha^d_jd_j-\sum_{i \in \mathcal{S}}\alpha^s_is_i-\sum_{\ell\in\ \mathcal{L}^{loc}}0\cdot f_\ell-\sum_{k\in\ \mathcal{K}}\sum_{\ell\in\ \mathcal{L}_k^{glo}}\alpha^f_{k} f_\ell-\sum_{t\in \mathcal{T}}\alpha^\xi_t\xi_t - \sum_{t\in \mathcal{T}}\alpha^y_t y_t\\
\label{eq:originalobj4} &=\sum_{j \in \mathcal{D}}\alpha^d_jd_j-\sum_{i \in \mathcal{S}}\alpha^s_is_i-\sum_{k\in\ \mathcal{K}}\sum_{\ell\in\ \mathcal{L}_k^{glo}}\alpha^{f}_{k} f_\ell-\sum_{t\in \mathcal{T}}\alpha^\xi_t\xi_t - \sum_{t\in \mathcal{T}}\alpha^y_t y_t\\
\label{eq:originalobj5} &=\sum_{j \in \mathcal{D}}\alpha^d_jd_j-\sum_{i \in \mathcal{S}}\alpha^s_is_i-\sum_{k\in\ \mathcal{K}}\alpha^{f}_{k} \hat{f}_k-\sum_{t\in \mathcal{T}}\alpha^\xi_t\xi_t - \sum_{t\in \mathcal{T}}\alpha^y_t y_t\\
\label{eq:originalobj6}&:=\overline{\varphi}.
\end{align}
\end{subequations}
The first expression \eqref{eq:originalobj1_1}  is the total welfare of the original supply chain problem. In  \eqref{eq:originalobj1_2} we define this in terms of local and global edges and in \eqref{eq:originalobj1_3} we partition the global edges.  In \eqref{eq:originalobj3} we eliminate the local flows (which is equivalent to setting a zero cost) and we define the cost of the global edges as $\alpha^f_{k}:=\min_{\ell\in\mathcal{L}_k^{glo}} \alpha^f_\ell$, which yields \eqref{eq:originalobj4}. We thus have that this expression provides an upper bound for the total welfare. In \eqref{eq:originalobj5} we rewrite this expression in terms of the aggregated global flows and we use this expression as the total welfare of the aggregated supply chain problem (which we define as $\overline{\varphi}$). By construction, this aggregation procedure yields $\varphi\leq \overline{\varphi}$ for any allocation $(s,d,f,\xi,y)$. 

The aggregation of the constraint set and the and objective give the approximate model:
\begin{subequations}
\begin{align}
\max_{(s,d,f,\xi,y)} & \; \sum_{j \in \mathcal{D}}\alpha^d_jd_j-\sum_{i \in \mathcal{S}}\alpha^s_is_i-\sum_{k\in\ \mathcal{K}}\alpha^f_{k} \hat{f}_k-\sum_{t\in \mathcal{T}}\alpha^\xi_t\xi_t - \sum_{t\in \mathcal{T}}\alpha^y_t y_t\\ 
\label{eq:clusterbalance}&\left(\sum_{i\in\mathcal{S}_{c,p}}s_i + \sum_{{k} \in
\mathcal{K}^{in}_{c, p}}\hat{f}_k \right) -\left(\sum_{j\in\mathcal{D}_{c,p}}d_j+\sum_{k\in \mathcal{K}^{out}_{c,p}}\hat{f}_k\right) + \sum_{t\in \mathcal{T}_{c}}\gamma_{t,p}\xi_t =0,\; (c,p)\in\mathcal{C}\times\mathcal{P}\\
&\; 0\leq d_j\leq \bar{d}_j,\; j\in \mathcal{D}\\
&\; 0\leq s_i\leq \bar{s}_i,\; i\in \mathcal{S}\\
&\label{eq:flow_hat_model2}\; 0\leq \hat{f}_k\leq \bar{\hat{f}}_{k},\ k\in\ \mathcal{K}\\
&\; 0\leq y_{t}\leq \bar{y}_{t},\; t\in \mathcal{T}\\
&\; 0\leq \xi_t\leq \bar{\xi}_t\cdot y_t,\; t\in \mathcal{T}
\end{align}
\end{subequations}
We denote the feasible region of this problem as $\bar{\mathcal{F}}$ and its optimal total welfare is $\bar{\varphi}^*$.  Because the associated constraint set is obtained by using linear combinations of the constraints of the feasible set $\mathcal{F}$, we have that $\mathcal{F}\subseteq \bar{\mathcal{F}}$. In other words, the feasible set of the aggregated  problem is a relaxation of the feasible set of the original problem. We also have that, by construction, $\varphi\leq \bar{\varphi}$ and thus $\varphi^*\leq \bar{\varphi}^*$ (the coarse problem provides a valid upper bound). 

\section{Numerical Case Study}

\begin{figure}[!htb]
\begin{center}
\includegraphics[height=6cm, width=16cm]{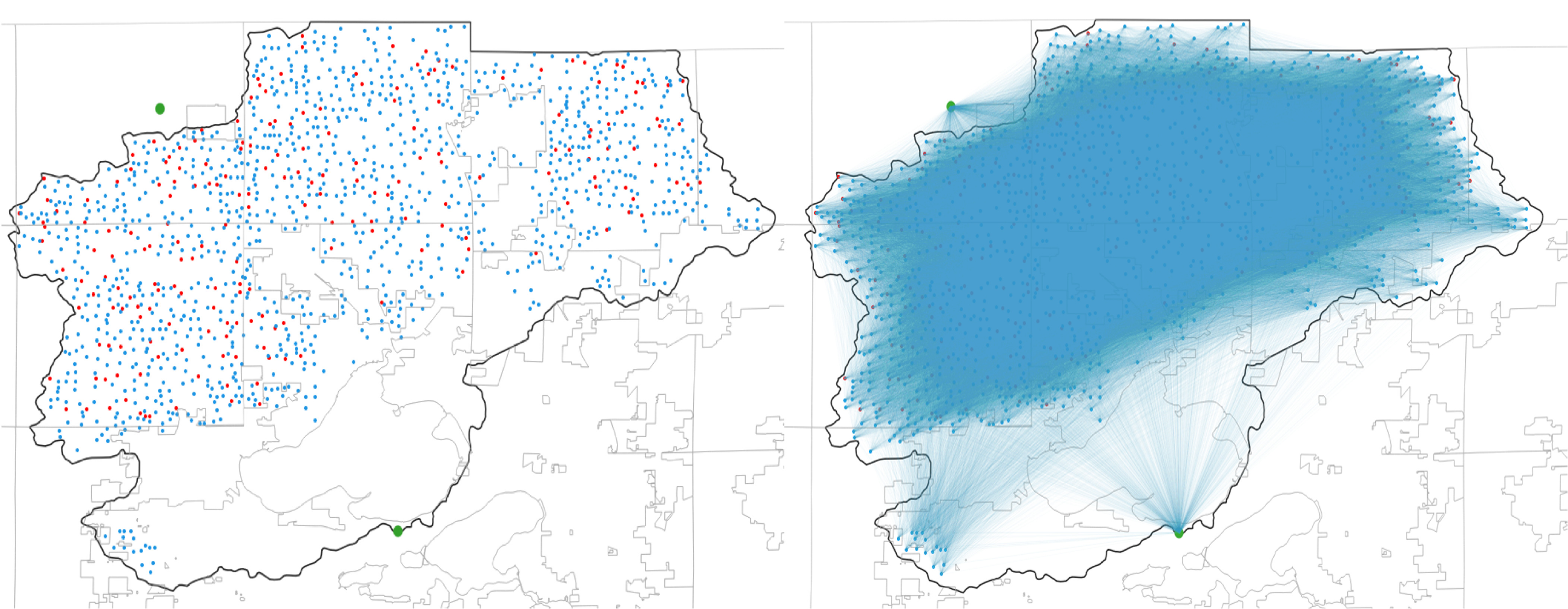}\caption{\mynote{\textcolor{black}{Left: spatial distribution of the supply chain nodes (dairy farms, croplands, and external companies). The blue dots include 1167 croplands, the red dots include 203 dairy farms, and the green dots are two consumers. Right: supply chain graph with 1.88 million edges}}}
\label{Fig.map}
\end{center}
\end{figure}

We demonstrate the proposed approach in a supply chain problem that arises in dairy waste management;  here, the goal is to deploy technologies to process the waste and obtain value-added projects. The study region is the upper Yahara watershed.  The study region contains 1,372 nodes, which correspond to 1,167 croplands, 203 dairy farms and a couple of external consumers (see Figure \ref{Fig.map}). The system includes a total of 20  products, which include the raw material (dairy manure) waste, intermediate products (e.g., solid waste and digestate) and final (value-added) products (e.g., struvite and pellets). A total of 12  technologies are considered for converting the raw material into intermediate and final products. We implement this  model in the algebraic modeling package JuMP and solve all problems using Gurobi 9.1.0. The code is executed on a computing server that contains a 17-core Intel(R) Xeon(R) CPU E5-2698 v3 processor running @2.30GHz.  All scripts and data needed to reproduce the results are available at: \url{https://github.com/zavalab/JuliaBox/tree/master/Graph_S%26C}. 

\begin{figure}[!ht]
     \begin{subfigure}{0.9\textwidth}
         \centering
         \includegraphics[width=\textwidth]{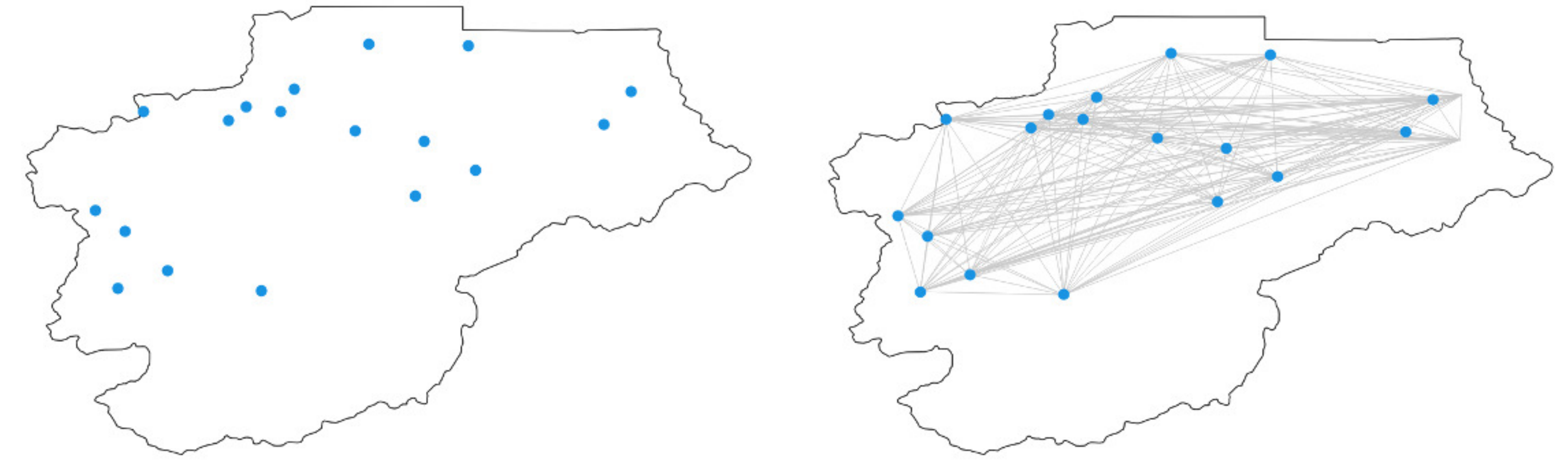} 
         \caption{Original graph with 20 nodes.}
         \label{Fig.agg_1.sub.1}
     \end{subfigure}
     \begin{subfigure}{0.9\textwidth}
         \centering
         \includegraphics[width=\textwidth]{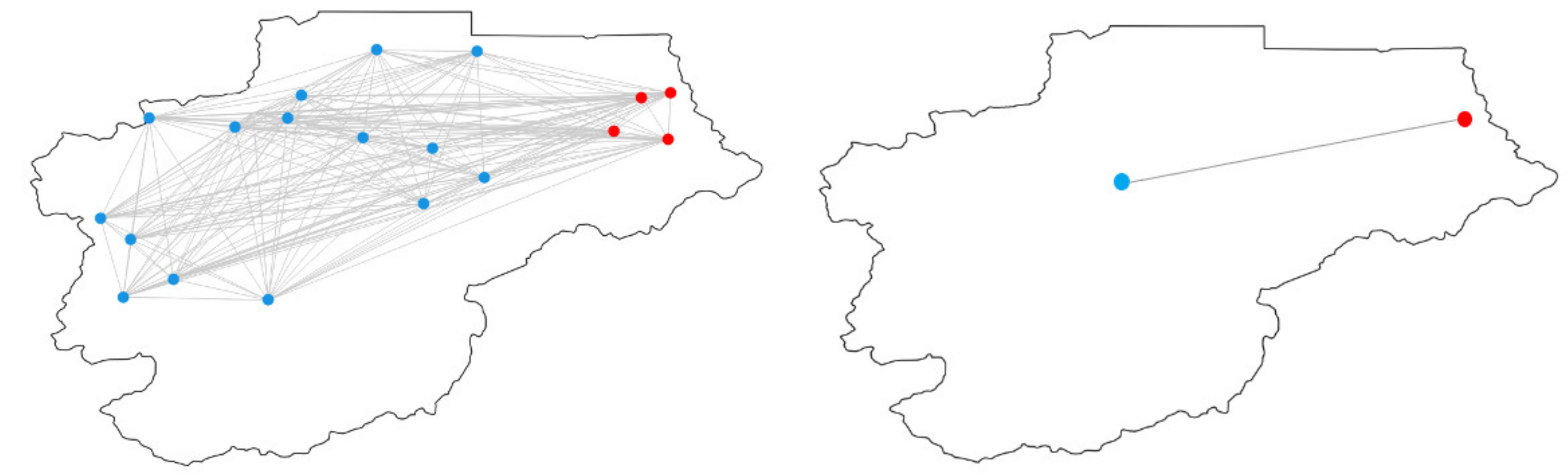}
         \caption{Coarsened graph with $C=2$ partitions and $K=1$ edges.}
         \label{Fig.agg_2.sub.2}
     \end{subfigure}
     
     \begin{subfigure}{0.9\textwidth}
         \centering
         \includegraphics[width=\textwidth]{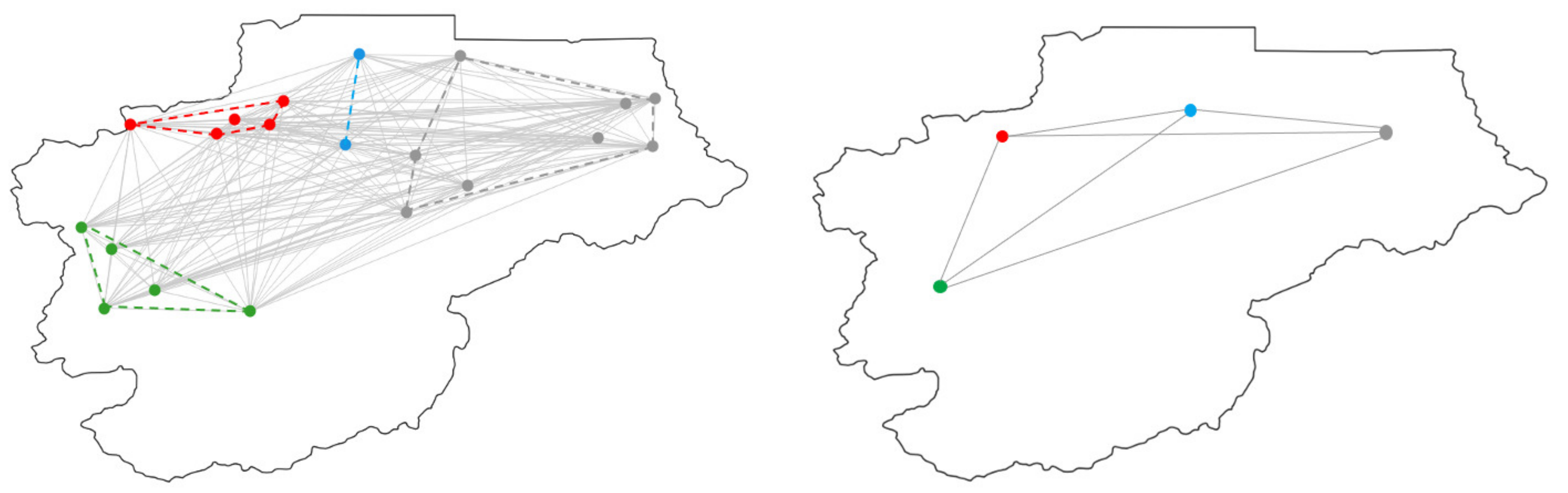}
         \caption{Graph with $C=4$ partitions and $K=6$ edges.}
         \label{Fig.agg_3.sub.3}
     \end{subfigure}
\caption{Graph coarsening for small case study. (a) Geographic distribution of nodes and the layout of the original network. (b) The 20 nodes are split into $C=2$ partitions and the partitions are aggregated to obtain 2 nodes (pivot nodes) and 1 edge. (c) The graph is split into $C=4$ partitions these are aggregated to obtain a graph with 4 nodes and 6 edges.}
\label{Fig.node_agg_case}
\end{figure}

We first consider a small-scale study to illustrate the concepts discussed; for this study, there are 20 nodes and $20^2=400$ edges in the supply chain graph. Figure \ref{Fig.node_agg_case} shows three levels of graph coarsening. Figure \ref{Fig.agg_1.sub.1} shows the spatial distribution of nodes and the layout of the original graph. In Figure \ref{Fig.agg_2.sub.2} we show the coarsening of the original supply chain to obtain $C=2$  partitions and $K=1$ edge. The local edges within each partition are eliminated and the global edges between two partitions are aggregated into a new edge between the partitions. In Figure \ref{Fig.agg_3.sub.3} we increase the number of partitions to $C=4$ and we obtain $K=6$ edges. Finally, at the third level, we do not perform any aggregation.

\begin{figure}[!htb]
\begin{centering}
     \begin{subfigure}{0.5\textwidth}
         \centering
         \includegraphics[width=\textwidth]{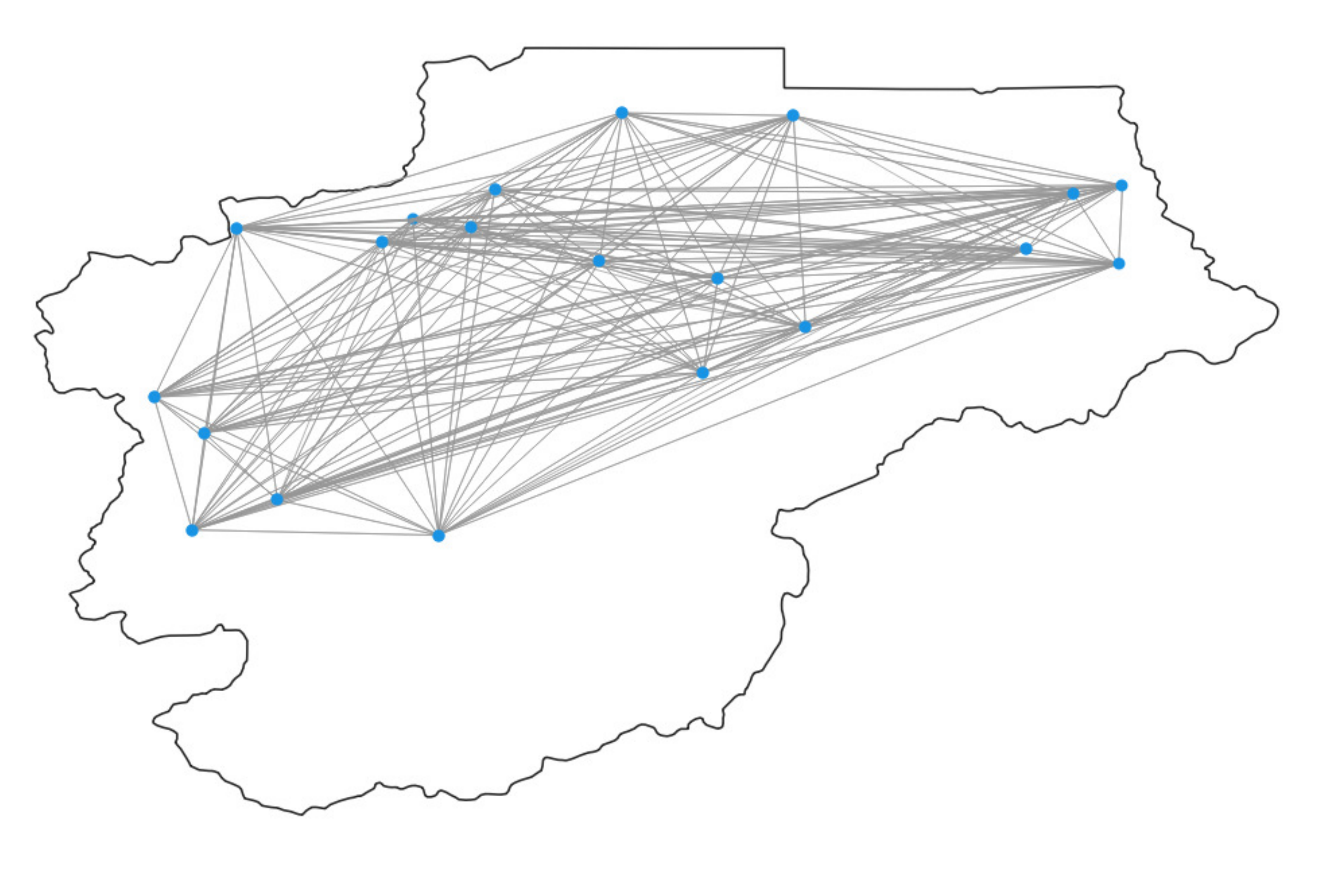}
         \caption{Original Graph with $|\mathcal{L}|=$ 400 edges.}
         \label{Fig.edge_case.sub.1}
     \end{subfigure}\\
     \begin{subfigure}{0.5\textwidth}
         \centering
         \includegraphics[width=\textwidth]{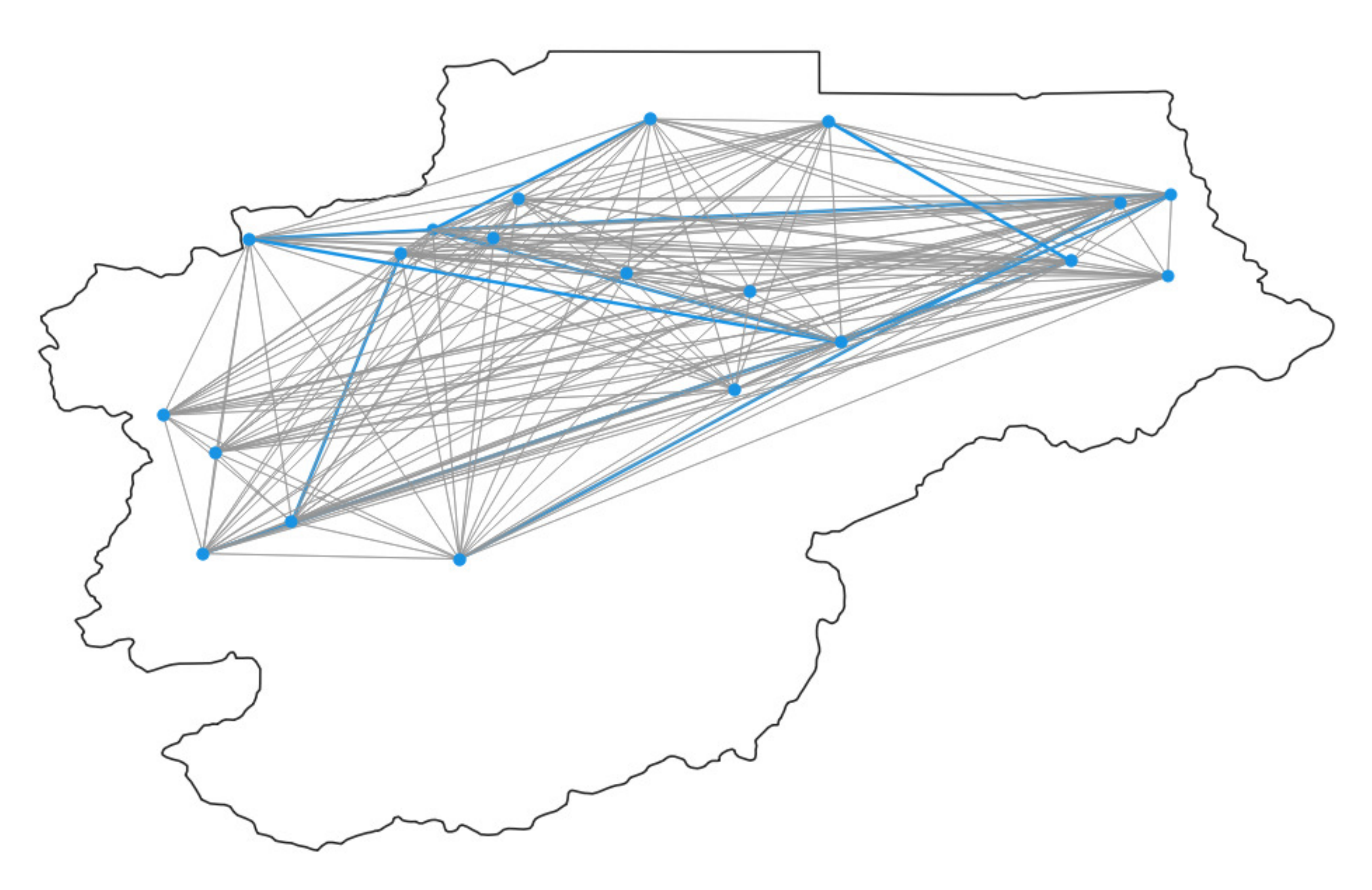}
         \caption{Sampled Graph with $|\mathcal{L}_a|=$10 active edges}
         \label{Fig.edge_case.sub.2}
     \end{subfigure}\\
     \begin{subfigure}{0.5\textwidth}
         \centering
         \includegraphics[width=\textwidth]{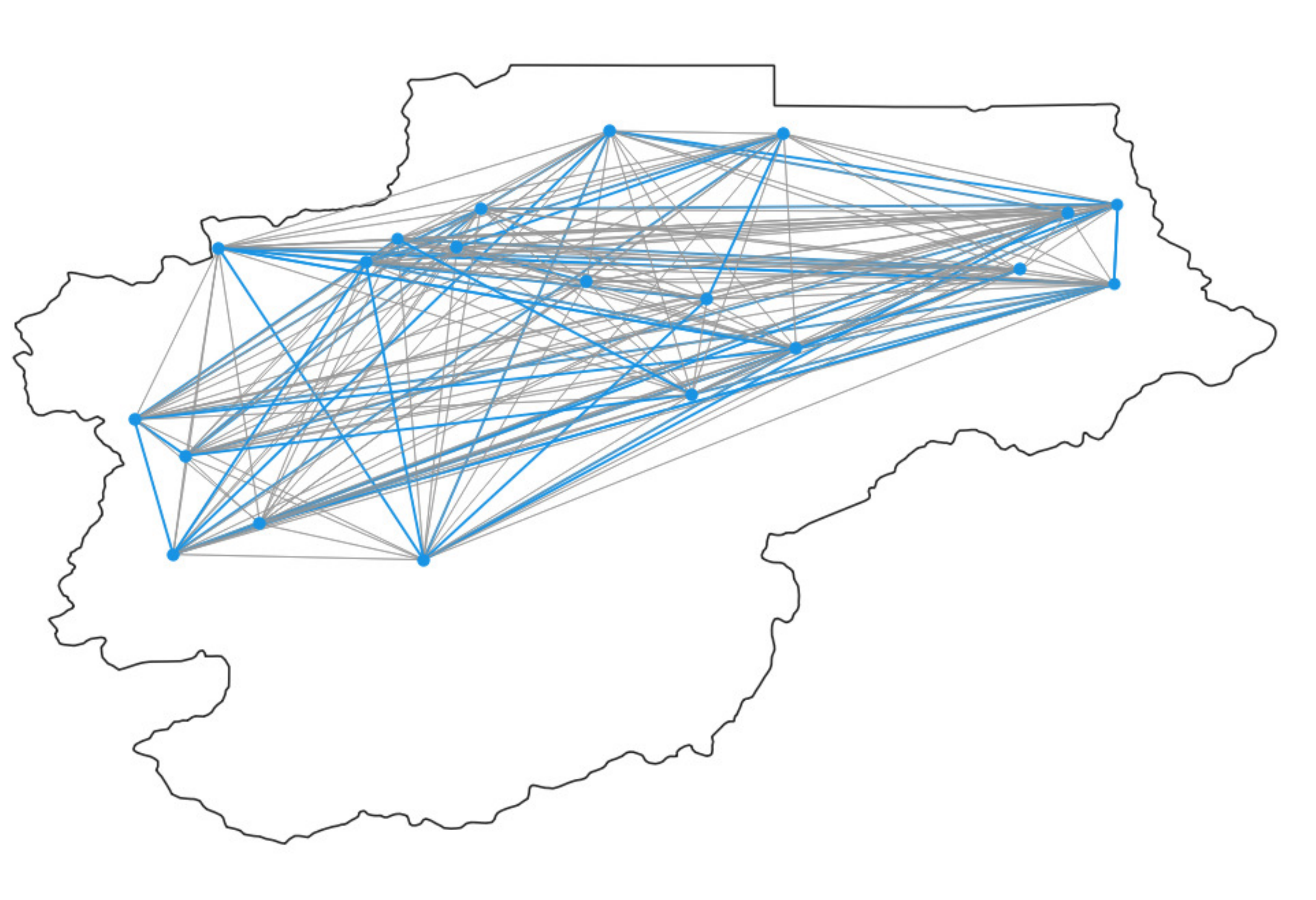}
         \caption{Sampled Graph with $|\mathcal{L}_a|=$40 active edges}
         \label{Fig.edge_case.sub.3}
     \end{subfigure}
\caption{Graph sampling for small case study. (a) Original graph with 20 nodes and 400 edges. (b) Sampled graph with 10 active edges. (c) Sampled graph with 40 active edges.}
\label{Fig.edge_case}
     \end{centering}
\end{figure}

Figure \ref{Fig.edge_case} illustrates the graph-sampling process. For the original graph, there are 400 edges (Figure \ref{Fig.edge_case.sub.1}). At the first level, $|\mathcal{L}_a|=$10 edges are randomly chosen to be active (Figure \ref{Fig.edge_case.sub.2}). We then increase the size of sample to $|\mathcal{L}_a|=$40 active edges (Figure \ref{Fig.edge_case.sub.3}). Here, we also show the remaining set of inactive edges. We can see that, as the number of sampled edges increases, the structure of the sampled graph resembles that of the original graph (even for a small number of edges). 

\begin{figure}[!htp]
\begin{center}
\includegraphics[scale=0.2]{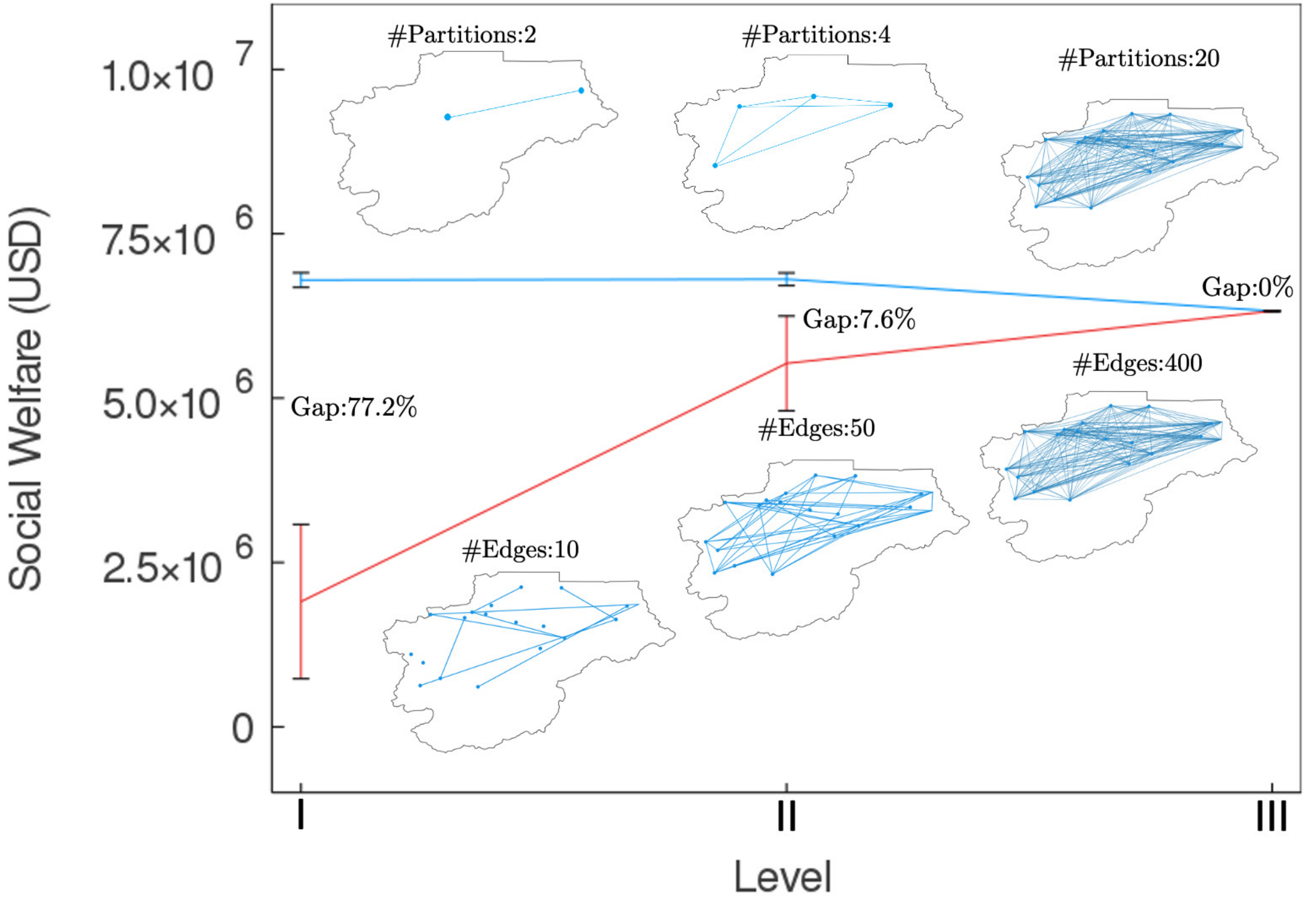}\caption{\mynote{\textcolor{black}{Computational results for small-scale study. Level I includes 10 sampled edges and 2 sampled nodes and achieves a 77.2\% optimality gap. Level II includes 40 sampled edges and 4 sampled nodes and achieves a 7.6\% gap. Level III includes 400 sampled edges and 20 sampled edges and achieves a 0\% gap.}}}
\label{Fig.result_smallcase}
\end{center}
\end{figure}

Figure \ref{Fig.result_smallcase} summarizes the computational results for the small-scale case study. We consider three levels of aggregation; for each level, we repeat the sampling process $|\Omega|=10$ times. Based on the samples, we compute 95\% confidence intervals for the lower and upper bounds. The first level uses 10 sampled edges (for graph sampling) and 2 sampled nodes (for coarsening). We recall that the number of sampled nodes equals the number of partitions and pivot nodes.  The approximate solution obtained in the first level achieves a 77.2\% optimality gap. Here, we can observe a large variability of the lower bound (graph sampling) but a small variability for the upper bound (coarsening). At the second level, we increase the number of sampled edges to 40 and the number of partitions to 4. We can see that this decreases the optimality gap by an order of magnitude (7.6\%) and nearly completely eliminates the variability of the lower bound. This confirms our intuition that, as the graph structure is better captured by the sampled graph, less variability in the total welfare that is observed. At the third level, we keep all the edges (400 sampled edges) and nodes in the graph (20 pivot nodes). As expected, as we increase the number of sampled edges and nodes, tighter bounds are obtained. This confirms that the gSC scheme is consistent and converges and that the variability collapses to zero. 

To illustrate the computational scalability of the gSC approach, we now apply it to a large-scale model that contains the 1,372 nodes, 1.8 million transport edges, 12 technologies, and 20 products. Table\ref{table1} summarizes the number of constraints and variables of the original model and of the approximate models obtained with graph-sampling and graph-coarsening. The  model contains more than 38 million continuous variables, 49,000 binary variables, and 78,000 constraints. We explored the performance of the proposed approach using three levels of approximation. At the third level of approximation, the graph-sampling approach used only 100,000 active edges (5.3\% of the total number of edges). This reduced the number of continuous variables by an order of magnitude (from 38 million to 2.4 million). The coarsening approach generates compact graphs by  eliminating the local edges and by collapsing global edges to single edges. This aggregation mechanism reduces the number of variables significantly; for instance, at the third level, the number of variables reduces from 38 million to just 500,000 (1.3\% of the variables). 
\begin{table}[!htb]
\centering
\caption{Size of the original, graph-sampling and graph-aggregation models}
\label{table1}
\begin{adjustbox}{max width=\textwidth}
\begin{tabular}{cccccccccccc}
\hline
                        & \multicolumn{5}{c}{Graph-Sampling}                                          &                         & \multicolumn{5}{c}{Graph-Aggregation}                                       \\ \cline{2-12} 
                        & \multicolumn{2}{c}{Variables} &           & \multicolumn{2}{c}{Constraints} &                         & \multicolumn{2}{c}{Variables} &           & \multicolumn{2}{c}{Constraints} \\ \cline{2-12} 
\multirow{-3}{*}{Level} & Continuous      & Binary      &           & Equality      & Inequality      &                         & Continuous      & Binary      &           & Equality      & Inequality      \\ \hline
1                       & 449,559         & 49,392      &           & 27,440        & 50,881          & {\color[HTML]{FFFFFF} } & 432,439         & 49,392      &           & 200           & 50,881          \\
2                       & 629,559         & 49,392      &           & 27,440        & 50,881          &                         & 439,239         & 49,392      &           & 400           & 50,881          \\
3                       & 2,429,559       & 49,392      & \textbf{} & 27,440        & 50,881          & \textbf{}               & 450,039         & 49,392      & \textbf{} & 600           & 50,881          \\ \hline
Original Model          & 38,093,703      & 49,392      &           & 27,440        & 50,881          &                         & 38,093,703      & 49,392      &           & 27,440        & 50,881          \\ \hline
\end{tabular}
\end{adjustbox}
\end{table}

\begin{figure}[!ht]
\begin{center}
\includegraphics[scale=0.2]{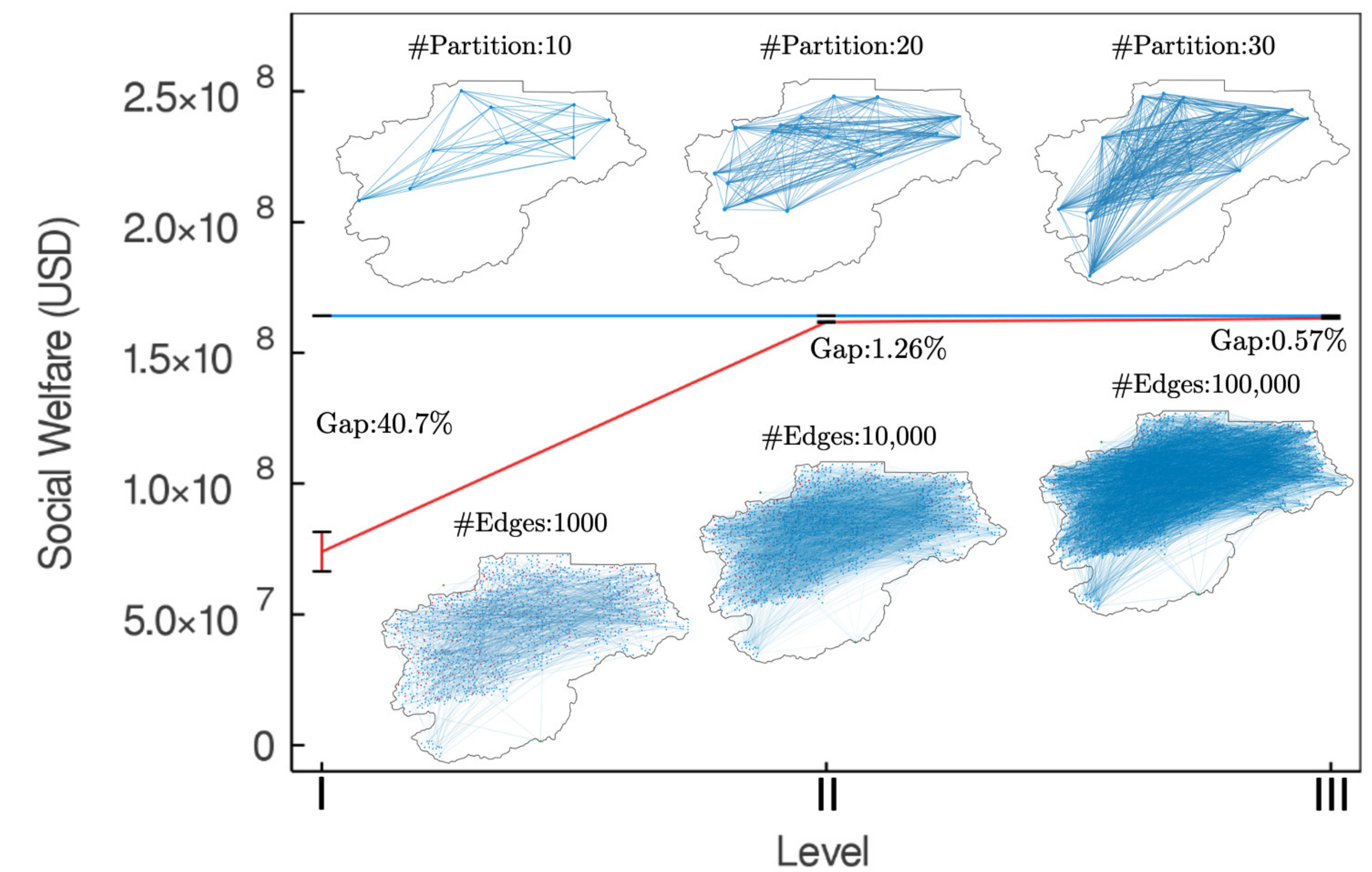}\caption{\mynote{\textcolor{black}{Computational results for large-scale study. Level I includes 1000 sampled edges and 10 sampled nodes and achieves a 40.7\% optimality gap. Level II includes 410,000 sampled edges and 20 sampled nodes and achieves a 1.26\% gap. Level III includes 100,000 sampled edges and 30 sampled edges and achieves a 0.57\% gap.}}}
\label{Fig.result_largecase}
\end{center}
\end{figure}

Table\ref{table2} shows the computational time and accuracy of the approximate solutions obtained. We can see that, by increasing the sample size, the standard deviation of the lower bounds decreases drastically. The computational time required for computing the lower bound in the third level is (on average) 735 sec. We can also see that tight upper bounds can be obtained after increasing the number of node samples from 10 to 30. The computational time required for solving the coarsened problem is (on average) 119 sec.  
Figure \ref{Fig.result_largecase} shows the lower and upper bounds of three levels. We again see that, by Increasing the sample size, the confidence intervals shrink and the bounds converge. Specifically, we see that the optimality gap decreases from 40.7\% (in the first level) to 0.57\% (in the third level). 
\\
\begin{table}[!htp]
\centering
\caption{Computational time and solution quality}
\label{table2}
\begin{adjustbox}{max width=\textwidth}
\begin{tabular}{cccccccccclc}
\hline
                        & \multicolumn{4}{c}{Graph-Sampling}         & \multicolumn{1}{l}{}    & \multicolumn{4}{c}{Graph-Aggregation}      &  &                       \\ \cline{2-5} \cline{7-10}
\multirow{-2}{*}{Level} & \#Edges  & Avg CPUs & Best LB   & Std of LB &                         & \#Partitions & Avg CPUs & Best UB   & Std of UB &  & \multirow{-2}{*}{Gap} \\ \hline
1                       & 1000    & 87       & \$97,366,828  & 13,139,210  &  & 10     & 60      & \$164,124,897 & 2445     &  & 40.7\%                \\
2                       & 10,000  & 101      & \$162,054,105  & 198,722   &                         & 20     & 79      & \$164,121,599 & 1623     &  & 1.26\%                 \\
3                       & 100,000 & 735      & \$163,179,868& 22,797     & \textbf{}               & 30     & 119     & \$164,114,871 & 1715    &  & 0.57\%                 \\ \hline
\end{tabular}
\end{adjustbox}
\end{table}

\begin{figure}[!htp]
\begin{center}
\includegraphics[height=8cm, width=11cm]{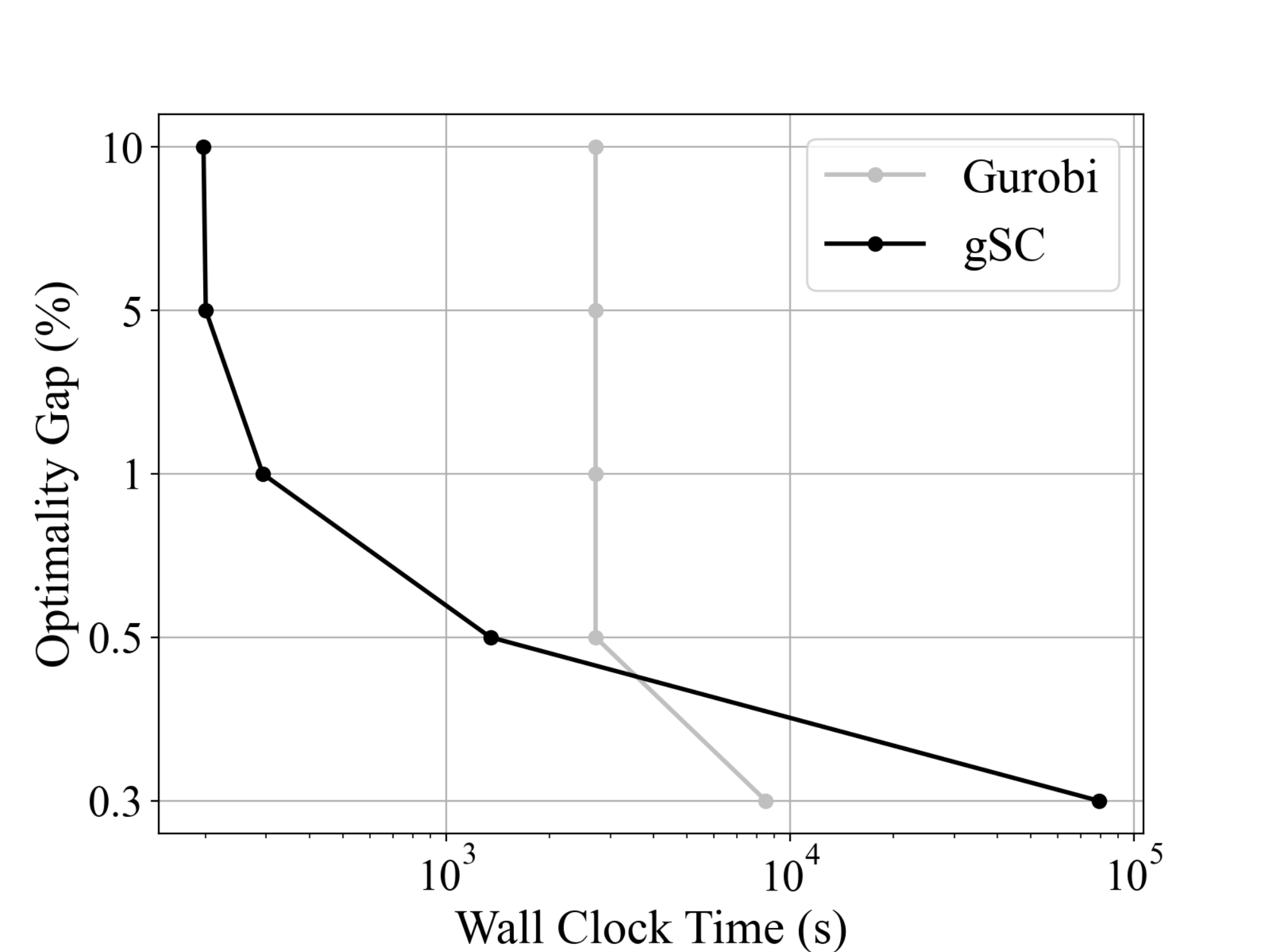}\caption{\mynote{\textcolor{black}{ Optimality gaps of gSC and Gurobi as a function of running time on the large-scale agricultural waste supply chain instance. To achieve 1\% optimality, gSC is ten times faster than Gurobi. To achieve 0.5\% of optimality, gSC is still two times faster. To achieve 0.3\% optimality, Gurobi is faster.  Neither Gurobi nor gSC can find any solution within 0.1\%  optimality within 24 hrs}}}
\label{Fig.time_gap}
\end{center}
\end{figure}

Figure \ref{Fig.time_gap}  \textcolor{black}{ shows the optimality gap as a function of wall-clock time.  The computational time of gSC includes the total time required for solving the problems in the three levels (lower\& upper bounds).  It takes more than 2000 sec for Gurobi to solve the relaxed linear program and takes almost 3000 sec to find a feasible solution.  On the other hand, gSC only requires 293 sec to find a feasible solution within a 1\% optimality gap.  Given the large size of the supply chain problem, we cannot find any solution within 0.1\%  optimality within 24 hrs by using either Gurobi or gSC. We also highlight that, in actual applications, one can also obtain a good estimate of the number of sampled edges and partitions used in gSC and thus do not have to account for the computational time of all levels.}\\

\begin{figure}[!htp]
\includegraphics[width=0.5\textwidth]{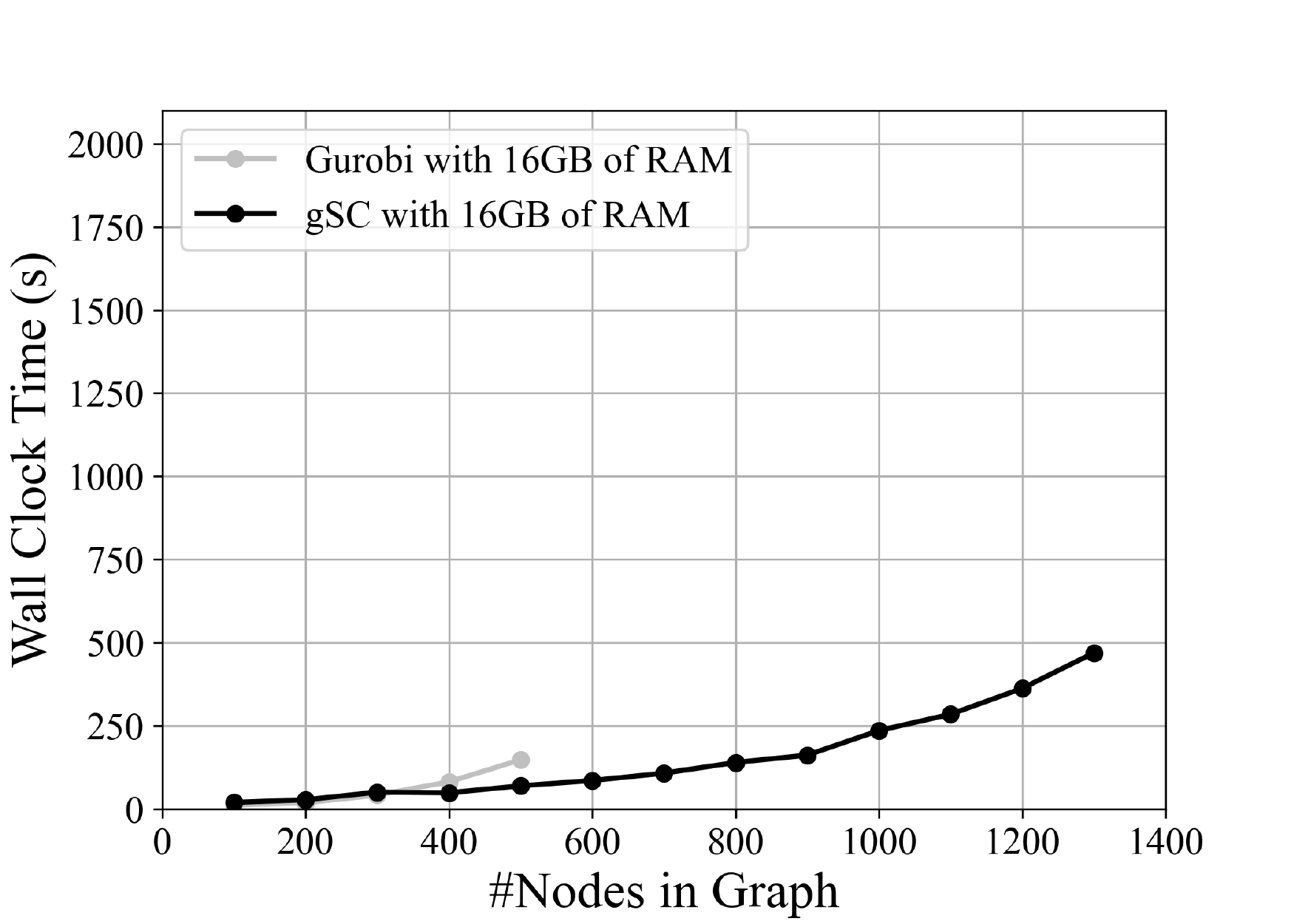}
\hspace{\fill}
\includegraphics[width=0.5\textwidth]{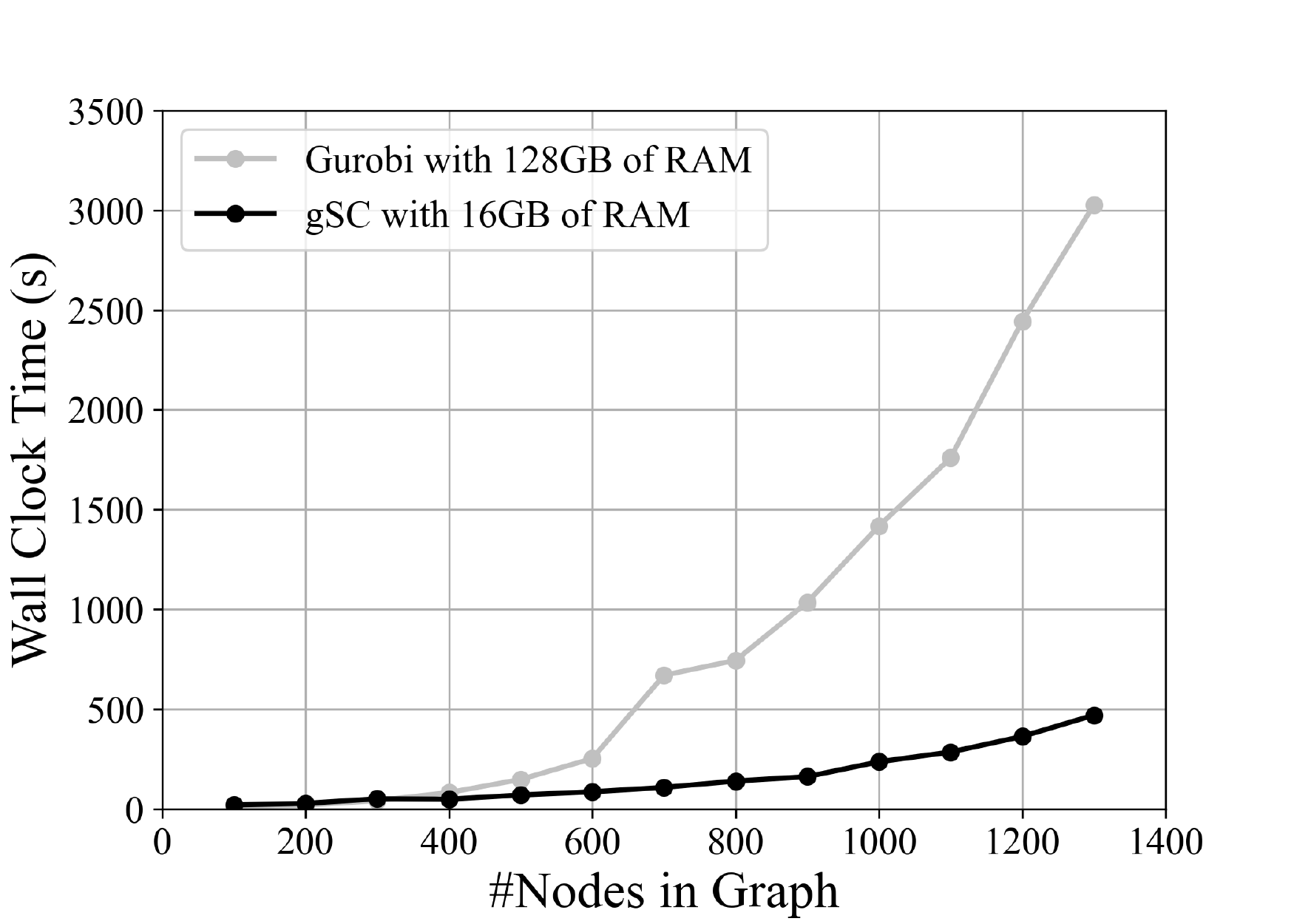}
\caption{Performance profile for test set containing 13 instances. Left: Gurobi and approximation scheme  with limited memory budget (16GB of RAM). Right: Using Gurobi with 128GB of RAM and using gSC algorithm with 16 GB of RAM. Gap of using Gurobi: <0.5\%, Gap of using gSC: <1\%}\label{fig:memory}
\end{figure}

We now compare the performance of the approximation schemes against Gurobi in models of different complexity. Here, we set the maximum RAM memory usage and test the largest model that Gurobi and gSC can handle with this memory budget. We then test Gurobi and gSC with a large amount of memory available and compare the computational time. Figure \ref{fig:memory} shows the computational times obtained with Gurobi and gSC for a tight memory budget of 16GB.  We start with a model with 100 nodes and progressively add nodes (this gives a test set with 13 problem instances). The largest instance in this test set is a model with 1,300 nodes. With a budget of 16 GB, only the first five instances can be solved by Gurobi. On the other hand, this memory budget is sufficient for finding approximate solutions with gSC in all instances. Specifically, the gSC finds high-quality solutions with <1\% gap. Our results also indicate that, in order to handle all 13 instances, Gurobi requires at least 128 GB of memory (which is about eight times higher than the memory usage of gSC). With sufficient memory allocated, Gurobi finds solutions with 0.5\% optimality gap for all 13 instances. As expected, the solution quality of the approximation scheme is lower than that achieved by Gurobi; however, the proposed approach can find approximate solutions much faster and/or with less memory. This type of flexibility is desirable, as it allows usage of different types of computing architectures. 

\section{Conclusions and Future Work}

We have presented a graph sampling and coarsening (gSC) scheme for computing high-quality approximate solutions for large-scale supply chain models. We show that this scheme can deliver lower bounds and upper bounds and can thus estimate the quality of the approximate solution. Our numerical results in a large model indicates that the proposed approach can tackle problems that are intractable to state-of-the-art solvers and can find good solutions under tighter memory and time budgets.  As part of future work, we are interested in exploring hierarchical coarsening schemes for handling large-scale supply chain models with space-time dynamics.

\section*{Acknowledgments}

We acknowledge funding from the U.S. Department of Agriculture, National Institute of Food and Agriculture, under grant number 2017-67003-26055.

\bibliography{fullpaper}

\begin{thebibliography}{10}

\bibitem{villa2002emerging}
Agostino Villa.
\newblock Emerging trends in large-scale supply chain management.
\newblock {\em International Journal of Production Research},
  40(15):3487--3498, 2002.

\bibitem{garcia2015supply}
Daniel~J Garcia and Fengqi You.
\newblock Supply chain design and optimization: Challenges and opportunities.
\newblock {\em Computers \& Chemical Engineering}, 81:153--170, 2015.

\bibitem{grossmann2012advances}
Ignacio~E Grossmann.
\newblock Advances in mathematical programming models for enterprise-wide
  optimization.
\newblock {\em Computers \& Chemical Engineering}, 47:2--18, 2012.

\bibitem{bok2000supply}
Jin-Kwang Bok, Ignacio~E Grossmann, and Sunwon Park.
\newblock Supply chain optimization in continuous flexible process networks.
\newblock {\em Industrial \& Engineering Chemistry Research}, 39(5):1279--1290,
  2000.

\bibitem{iyer1998bilevel}
Ramaswamy~R Iyer and Ignacio~E Grossmann.
\newblock A bilevel decomposition algorithm for long-range planning of process
  networks.
\newblock {\em Industrial \& Engineering Chemistry Research}, 37(2):474--481,
  1998.

\bibitem{dogan2006decomposition}
Muge~Erdirik Dogan and Ignacio~E Grossmann.
\newblock A decomposition method for the simultaneous planning and scheduling
  of single-stage continuous multiproduct plants.
\newblock {\em Industrial \& engineering chemistry research}, 45(1):299--315,
  2006.

\bibitem{pishvaee2014accelerated}
Mir~Saman Pishvaee, Jafar Razmi, and Seyed~Ali Torabi.
\newblock An accelerated benders decomposition algorithm for sustainable supply
  chain network design under uncertainty: A case study of medical needle and
  syringe supply chain.
\newblock {\em Transportation Research Part E: Logistics and Transportation
  Review}, 67:14--38, 2014.

\bibitem{oliveira2014accelerating}
Fabricio Oliveira, Ignacio~E Grossmann, and Silvio Hamacher.
\newblock Accelerating benders stochastic decomposition for the optimization
  under uncertainty of the petroleum product supply chain.
\newblock {\em Computers \& Operations Research}, 49:47--58, 2014.

\bibitem{uster2007benders}
Halit {\"U}ster, Gopalakrishnan Easwaran, Elif Ak{\c{c}}ali, and Sila
  {\c{C}}etinkaya.
\newblock Benders decomposition with alternative multiple cuts for a
  multi-product closed-loop supply chain network design model.
\newblock {\em Naval research logistics (NRL)}, 54(8):890--907, 2007.

\bibitem{santoso2005stochastic}
Tjendera Santoso, Shabbir Ahmed, Marc Goetschalckx, and Alexander Shapiro.
\newblock A stochastic programming approach for supply chain network design
  under uncertainty.
\newblock {\em European Journal of Operational Research}, 167(1):96--115, 2005.

\bibitem{jackson2003temporal}
Jennifer~R Jackson and Ignacio~E Grossmann.
\newblock Temporal decomposition scheme for nonlinear multisite production
  planning and distribution models.
\newblock {\em Industrial \& engineering chemistry research},
  42(13):3045--3055, 2003.

\bibitem{oliveira2013lagrangean}
Fabricio Oliveira, Vijay Gupta, Silvio Hamacher, and Ignacio~E Grossmann.
\newblock A lagrangean decomposition approach for oil supply chain investment
  planning under uncertainty with risk considerations.
\newblock {\em Computers \& Chemical Engineering}, 50:184--195, 2013.

\bibitem{terrazas2011temporal}
Sebastian Terrazas-Moreno, Philipp~A Trotter, and Ignacio~E Grossmann.
\newblock Temporal and spatial lagrangean decompositions in multi-site,
  multi-period production planning problems with sequence-dependent
  changeovers.
\newblock {\em Computers \& Chemical Engineering}, 35(12):2913--2928, 2011.

\bibitem{terrazas2011multiscale}
Sebastian Terrazas-Moreno and Ignacio~E Grossmann.
\newblock A multiscale decomposition method for the optimal planning and
  scheduling of multi-site continuous multiproduct plants.
\newblock {\em Chemical Engineering Science}, 66(19):4307--4318, 2011.

\bibitem{sousa2011global}
Rui~T Sousa, Songsong Liu, Lazaros~G Papageorgiou, and Nilay Shah.
\newblock Global supply chain planning for pharmaceuticals.
\newblock {\em chemical engineering research and design}, 89(11):2396--2409,
  2011.

\bibitem{van2001lagrangean}
Susara~A van~den Heever, Ignacio~E Grossmann, Sriram Vasantharajan, and
  Krisanne Edwards.
\newblock A lagrangean decomposition heuristic for the design and planning of
  offshore hydrocarbon field infrastructures with complex economic objectives.
\newblock {\em Industrial \& engineering chemistry research},
  40(13):2857--2875, 2001.

\bibitem{barmann2015solving}
Andreas B{\"a}rmann, Frauke Liers, Alexander Martin, Maximilian Merkert,
  Christoph Thurner, and Dieter Weninger.
\newblock Solving network design problems via iterative aggregation.
\newblock {\em Mathematical Programming Computation}, 7(2):189--217, 2015.

\bibitem{chen2017large}
Annie I-An Chen.
\newblock {\em Large-scale optimization in online-retail inventory management}.
\newblock PhD thesis, Massachusetts Institute of Technology, 2017.

\bibitem{wu2021predictive}
Yaqing Wu, Christos~T Maravelias, Michael~J Wenzel, Mohammad~N ElBsat, and
  Robert~T Turney.
\newblock Predictive maintenance scheduling optimization of building heating,
  ventilation, and air conditioning systems.
\newblock {\em Energy and Buildings}, 231:110487, 2021.

\bibitem{zipkin1982exact}
Paul Zipkin.
\newblock Exact and approximate cost functions for product aggregates.
\newblock {\em Management Science}, 28(9):1002--1012, 1982.

\bibitem{zipkin1982transportation}
Paul Zipkin.
\newblock Transportation problems with aggregated destinations when demands are
  uncertain.
\newblock {\em Naval Research Logistics Quarterly}, 29(2):257--270, 1982.

\bibitem{clautiaux2017iterative}
Fran{\c{c}}ois Clautiaux, Said Hanafi, Rita Macedo, Marie-Emilie Voge, and
  Cl{\'a}udio Alves.
\newblock Iterative aggregation and disaggregation algorithm for
  pseudo-polynomial network flow models with side constraints.
\newblock {\em European Journal of Operational Research}, 258(2):467--477,
  2017.

\bibitem{shetty1987solving}
CM~Shetty and Richard~W Taylor.
\newblock Solving large-scale linear programs by aggregation.
\newblock {\em Computers \& Operations Research}, 14(5):385--393, 1987.

\bibitem{zipkin1980variable}
Paul~H Zipkin.
\newblock Bounds on the effect of aggregating variables in linear programs.
\newblock {\em Operations Research}, 28(2):403--418, 1980.

\bibitem{litvinchev2013aggregation}
Igor Litvinchev and Vladimir Tsurkov.
\newblock {\em Aggregation in large-scale optimization}, volume~83.
\newblock Springer Science \& Business Media, 2013.

\bibitem{zipkin1980nodes}
Paul~H Zipkin.
\newblock Bounds for aggregating nodes in network problems.
\newblock {\em Mathematical programming}, 19(1):155--177, 1980.

\bibitem{verweij2003sample}
Bram Verweij, Shabbir Ahmed, Anton~J Kleywegt, George Nemhauser, and Alexander
  Shapiro.
\newblock The sample average approximation method applied to stochastic routing
  problems: a computational study.
\newblock {\em Computational optimization and applications}, 24(2):289--333,
  2003.

\bibitem{kleywegt2002sample}
Anton~J Kleywegt, Alexander Shapiro, and Tito Homem-de Mello.
\newblock The sample average approximation method for stochastic discrete
  optimization.
\newblock {\em SIAM Journal on Optimization}, 12(2):479--502, 2002.

\bibitem{schutz2009supply}
Peter Sch{\"u}tz, Asgeir Tomasgard, and Shabbir Ahmed.
\newblock Supply chain design under uncertainty using sample average
  approximation and dual decomposition.
\newblock {\em European journal of operational research}, 199(2):409--419,
  2009.

\bibitem{bidhandi2017accelerated}
Hadi~Mohammadi Bidhandi and Jonathan Patrick.
\newblock Accelerated sample average approximation method for two-stage
  stochastic programming with binary first-stage variables.
\newblock {\em Applied Mathematical Modelling}, 41:582--595, 2017.

\bibitem{li2018sample}
Xueping Li and Kaike Zhang.
\newblock A sample average approximation approach for supply chain network
  design with facility disruptions.
\newblock {\em Computers \& Industrial Engineering}, 126:243--251, 2018.

\bibitem{sampat2019coordinated}
Apoorva~M Sampat, Yicheng Hu, Mahmoud Sharara, Horacio Aguirre-Villegas,
  Gerardo Ruiz-Mercado, Rebecca~A Larson, and Victor~M Zavala.
\newblock Coordinated management of organic waste and derived products.
\newblock {\em Computers \& chemical engineering}, 128:352--363, 2019.

\bibitem{tominac2021economic}
Philip~A Tominac and Victor~M Zavala.
\newblock Economic properties of multi-product supply chains.
\newblock {\em Computers \& Chemical Engineering}, 145:107157, 2021.

\end{thebibliography}
\end{document}